\documentclass[11pt]{amsart}

\usepackage{amsmath,amssymb,amscd}
\usepackage{enumerate}
\usepackage{xcolor}

\usepackage{amsmath,amssymb}

\usepackage{amsxtra, amsmath}
\usepackage{amssymb, amscd}
\usepackage{graphicx}
\usepackage{mathrsfs}

\usepackage{mathtools}

\newtheorem{thm}{Theorem}[section]
\newtheorem{prop}[thm]{Proposition}

\newtheorem{cor}[thm]{Corollary}

\theoremstyle{definition}
\newtheorem{definition}[thm]{Definition}

\theoremstyle{remark}
\newtheorem{remark}[thm]{Remark}
\newtheorem{example}[thm]{Example}

\usepackage{dsfont}
\allowdisplaybreaks

\title{ The algebra  $\mathcal{D}(W)$ via  strong Darboux transformations}
\author{Ignacio Bono Parisi}
\author{Ines Pacharoni}
\thanks{This paper was partially supported by SeCyT-UNC, CONICET, PIP 1220150100356.}

\address{CIEM-FaMAF\\ Universidad Nacional de C\'or\-do\-ba\\
CP 5000, C\'or\-do\-ba,  Argentina}
\email{ignacio.bono@unc.edu.ar, ines.pacharoni@unc.edu.ar}

\begin{document}
\begin{abstract}
The Matrix Bochner Problem aims to classify weight matrices $W$ such that the algebra $\mathcal D(W)$, of all differential operators that have a sequence of matrix-valued orthogonal polynomials for $W$ as eigenfunctions, contains a second-order differential operator. 
In \cite{CY18} it is proven that, under certain assumptions, the solutions to the Matrix Bochner Problem can be obtained through a noncommutative bispectral Darboux transformation of some classical scalar weights.

The main aim of this paper is to introduce the concept of strong Darboux transformation among weight matrices and explore the relationship between the algebras $\mathcal{D}(W)$ and $\mathcal{D}(\widetilde{W})$ when $\widetilde{W}$ is a strong Darboux transformation of $W$.  Starting from a direct sum of classical scalar weights $\widetilde W$, and leveraging our complete knowledge of the algebra of $\mathcal D(\widetilde W)$, we can easily 
determine the algebra $\mathcal D(W)$ of a weight $W$ that is a strong Darboux transformation of $\widetilde W$.
\end{abstract} 

\maketitle

\section{Introduction}
The theory of matrix-valued orthogonal polynomials starts with the work of Krein, see \cite{K49, K71}. If one is considering possible applications of these matrix polynomials, it is natural to concentrate on those cases where some extra property holds. A. Durán, in \cite{D97}, posed the problem of characterizing those positive definite matrix-valued weights whose matrix-valued orthogonal polynomials will be common eigenfunctions of some symmetric second-order differential operator with matrix coefficients. 
This was an important generalization to the matrix-valued case of the problem originally considered by S. Bochner in 1929, see \cite{B29}.
S. Bochner proved that, up to an affine change of coordinates, the only scalar weights satisfying these properties are the classical Hermite, Laguerre, and Jacobi weights: $w(x)=e^{-x^2}$, $w_{\alpha}(x)=x^\alpha e^{-x}$, and $w_{\alpha,\beta}(x)=(1-x)^\alpha(1 + x)^\beta$ respectively.

The problem of finding weight matrices $W(x)$ of size $N \times N$  such that the associated sequence of orthogonal matrix polynomials are eigenfunctions of a second-order differential operator is known as the Matrix Bochner Problem. 
  In the matrix-valued case,  the first nontrivial solutions of this problem were given, in   \cite{GPT01},  \cite{GPT02}, \cite{G03}, \cite{DG04}, mainly as a consequence of the study of matrix-valued spherical functions. 
In the past twenty years, a large amount of examples have been found. 
See \cite{GPT03a}, \cite{GPT03b}, \cite{GPT05}, \cite{CG06}, \cite{PT07}, \cite{DG07},\cite{P08}, \cite{PR08}, \cite{PZ16}, \cite{KPR12}. 
However, the explicit construction of such examples is, even nowadays, not an easy task.

After the appearance of the first examples, some works focused on studying the algebra $\mathcal{D}(W)$ of all differential operators that have a sequence of matrix-valued orthogonal polynomials for $W$ as eigenfunctions. 
This algebra is studied both from a general standpoint and for specific weights $W$, see \cite{GT07}, \cite{CG06}, \cite{PZ16},\cite{T11}, \cite{Z15}. 
The Matrix Bochner Problem can be read as to find the weight matrices $W$ such that its algebra $\mathcal D(W)$ contains a second-order differential operator.

 W. R. Casper and M. Yakimov in \cite{CY18}, studying this algebra $\mathcal D(W)$,  made a significant breakthrough towards solving the Matrix Bochner Problem. They proved that, under certain assumptions on the algebra $\mathcal{D}(W)$, the solutions to the Matrix Bochner Problem can be obtained through a noncommutative bispectral Darboux transformation of a direct sum of classical scalar weights. 
 Although this approach does not provide a complete solution to the Matrix Bochner Problem, as proven in \cite{BP23} and \cite{BP24}, it encourages us to further explore the connection between the algebra $\mathcal{D}(W)$ and Darboux transformations among weights.

Starting from a direct sum of classical scalar weights, one can use a Darboux transformation to obtain new irreducible weights whose algebras will be non-trivial, as seen in \cite{C18}.  
In \cite{BP23b} we study the Darboux transformations between direct sums of classical scalar weights of Hermite, Laguerre, and Jacobi type and we use these results to determine the algebra $\mathcal D(W)$ when $W$ is a direct sum of classical scalar weights of the same type.  

Calculating the algebra of an irreducible weight of size $N$ is a highly complex problem. In the literature, there are several works where authors consider an irreducible weight $W$, and by using the equations of symmetry for an operator $D$, they calculate operators up to a certain order and then conjecture possible generators of the algebra. 

\smallskip

We say that $\widetilde{W}$ is a strong Darboux transformation of $W$ if there exists a differential operator $D \in \mathcal{D}(W)$  that can be factorized as $D = \mathcal{V}\mathcal{N}$, with $\mathcal{N} = \widetilde{W}\mathcal{V}^{\ast}W^{-1}$, 
$\mathcal{V}$ a degree-preserving differential operator with nonsingular leading coefficient, and such that 
$$\langle P\cdot \mathcal{V}, Q\rangle_{\widetilde{W}} = \langle P , Q\cdot \mathcal{N} \rangle_{W},  \quad \text{ for all }  
P,Q \in \operatorname{Mat}_{N}(\mathbb{C}[x]).$$ (See Definition \ref{strong darb}). 
 
The specific factorization of a differential operator $D\in \mathcal D(W)$, in the form $  D=  \mathcal{V}  \widetilde{W}\mathcal{V}^{\ast}W^{-1}$ allows us to gather a significant amount of information about $\widetilde{W}$ by knowing that of $W$. 
For instance, 
we  see that 
$\mathcal{V}$ maps the sequence of monic orthogonal polynomials of $W$ to a sequence of orthogonal polynomials for $\widetilde{W}$.
 (See  Proposition \ref{strong1}). Also
their algebras are strongly related: for any differential operator $\mathcal{A} \in \mathcal{D}(W)$, we have $\mathcal{B} = \mathcal{N}\mathcal{A}\mathcal{V} \in \mathcal{D}(\widetilde{W})$ and  moreover, if $\mathcal{A}$ is $W$-symmetric, then $\mathcal{B}$ is $\widetilde{W}$-symmetric.   
In particular, the order of $\mathcal{B}$ has the same parity as the order of $\mathcal{A}$.  
(See Theorems \ref{main} and  \ref{parity}). 

To  illustrate this deep connection between the algebras $\mathcal D(W)$ and $\mathcal D(\widetilde W)$ when $W$ and $
\widetilde W$ are strong Darboux transformations of each other, we consider several examples of irreducible weights $W$. We prove that they are a strong Darboux transformation of a direct sum of scalar weights, and thus we determine the algebra $\mathcal D(W)$ in a very simple way.    
 
In Section \ref{Laguerre-Sect}, 
we consider a Laguerre-type irreducible weight $W(x) = e^{-x}x^{\alpha} \begin{psmallmatrix} x(1+a^{2}x) && ax \\ax && 1 \end{psmallmatrix}$.  
This example first appeared in \cite{DI08}. Based on computational calculations, the authors conjectured that certain operators were generators of the algebra $D(W)$ and provided some algebraic relations for these generators. 
We prove in Proposition \ref{ddd} that $W$ is a strong Darboux transformation of the direct sum of classical scalar Laguerre weights $\begin{psmallmatrix}e^{-x}x^{\alpha+1} && 0 \\ 0 && e^{-x}x^{\alpha} \end{psmallmatrix}$. This allows us to determine the algebra $\mathcal{D}(W)$ and a sequence of orthogonal polynomials for $W$.

In \cite{DG07}, the authors found that for size $2\times 2$, there are only 
three families of irreducible weight matrices of Hermite-type such that their algebras admit a second-order differential operator with leading coefficient $F_{2}(x) = I$. In Section \ref{Hermite Sect}, we study these three families. 
For  $W_{1}(x) = e^{-x^{2}} \begin{psmallmatrix} 1 +a^{2}x^{4} && ax^{2} \\ ax^{2} && 1 \end{psmallmatrix}$, there is no information or conjectures available in the literature regarding the algebra $\mathcal D(W_1)$. We prove that $W_{1}$ is a strong Darboux transformation of a direct sum of classical scalar Hermite weights. Thus we determine the algebra $\mathcal{D}(W_{1})$ and a sequence of orthogonal polynomials. Additionally, we present some algebraic relations for its generators, the three-term recurrence relation of its sequence of orthogonal polynomials, and we calculate the generators of the center $\mathcal{Z}(W_{1})$ of $\mathcal{D}(W_{1})$.
The second family is $W_{2}(x) = e^{-x^{2}}\begin{psmallmatrix} 1 + a^{2}x^{2} && ax \\ ax && 1 \end{psmallmatrix}$. The algebra $\mathcal D(W_2) $ was already determined by J. Tirao in \cite{T11}. We prove that $W_{2}$ is a strong Darboux transformation of $w(x) = e^{-x^{2}}I$, allowing us to provide a shorter argument for determining its algebra $\mathcal{D}(W_{2})$. Finally, the third family   is $W_{3}(x) = e^{-x^{2}} \begin{psmallmatrix} e^{2bx} + a^{2}x^{2} && ax \\ ax && 1 \end{psmallmatrix}$. It was considered in \cite{BP23}, where we establish that $W_{3}$ is not a Darboux transformation of any direct sum of classical scalar weights. We recover this result using the relationship that would exist between the algebras if this weight were a Darboux transformation of a direct sum of scalar weights.

In the final section, we present two novel families of Hermite-type weight matrices which are a strong Darboux transformation of a direct sum of classical scalar Hermite weights. It is noteworthy that these weight matrices while being Darboux transformations of classical weights, are not solutions to the Matrix Bochner Problem. The algebra of the first example contains a fourth-order differential operator and contains no differential operators of order less than $4$. Similarly, the algebra of the second example contains a sixth-order differential operator and contains no differential operators of order less than $6$. 

\section{Matrix valued orthogonal polynomials and the algebra $\mathcal D(W)$}\label{sec-MOP}

 Let $W=W(x)$ be a weight matrix of size $N$ on the real line, that is, a complex $N\times N$ matrix-valued smooth function on the interval $(x_0,x_1)$ such that $W(x)$ is positive definite almost everywhere and with finite moments of all orders. Let $\operatorname{Mat}_N(\mathbb{C})$ be the algebra of all $N\times N$ complex matrices and let $\operatorname{Mat}_N(\mathbb{C}[x])$ be the algebra of polynomials in the indeterminate $x$ with coefficients in $\operatorname{Mat}_N(\mathbb{C})$. We consider the following Hermitian sesquilinear form in the linear space $\operatorname{Mat}_N(\mathbb{C}[x])$
\begin{equation*}
  \langle P,Q \rangle =  \langle P,Q \rangle_W = \int_{x_0}^{x_1} P(x) W(x) Q(x)^*\,dx.
\end{equation*}

Given a weight matrix $W$ one can construct sequences 
$\{Q_n\}_{n\in\mathbb{N}_0}$ of matrix-valued orthogonal polynomials, i.e. the 
$Q_n$ are polynomials of degree $n$ with nonsingular leading coefficient and $\langle Q_n,Q_m\rangle=0$ for $n\neq m$.
We observe that there exists a unique sequence of monic orthogonal polynomials $\{P_n\}_{n\in\mathbb{N}_0}$ in $\operatorname{Mat}_N(\mathbb{C}[x])$.
By following a standard argument (see \cite{K49} or \cite{K71}) one shows that the monic orthogonal polynomials $\{P_n\}_{n\in\mathbb{N}_0}$ satisfy a three-term recursion relation
\begin{equation}\label{ttrr}
    x P_n(x)=P_{n+1}(x) + B_{n}P_{n}(x)+ C_nP_{n-1}(x), \qquad n\in\mathbb{N}_0,
\end{equation}
where $P_{-1}=0$ and $B_n, C_n$ are matrices depending on $n$ and not on $x$.

Along this paper, we consider that an arbitrary matrix differential operator
\begin{equation}\label{D2}
  {D}=\sum_{i=0}^s \partial ^i F_i(x),\qquad \partial=\frac{d}{dx},
\end{equation}
acts on the right   on  a matrix-valued function $P $ i.e.
$(P{D})(x)=\sum_{i=0}^s \partial ^i (P)(x)F_i(x).$

We consider the algebra of these operators with polynomial coefficients $$\operatorname{Mat}_{N}(\Omega[x])=\Big\{D = \sum_{j=0}^{n} \partial^{j}F_{j}(x) \, : F_{j} \in \operatorname{Mat}_{N}(\mathbb{C}[x]) \Big \}.$$
\noindent 
More generally, when necessary, we will also consider $\operatorname{Mat}_{N}(\Omega[[x]])$, the set of all differential operators with coefficients in  $\mathbb{C}[[x]]$, the ring of power series with coefficients in $\mathbb{C}$.

\begin{prop}[\cite{GT07}, Propositions 2.6 and 2.7]\label{eigenvalue-prop}
  Let $W=W(x)$ be a weight matrix of size $N$ and let $\{P_n\}_{n\geq 0}$ be  the sequence of monic orthogonal polynomials in $\operatorname{Mat}_N(\mathbb{C}[x])$. If $D$ is differential operator of order $s$, as in \eqref{D2}, such that
  $$P_nD=\Lambda_n P_n, \qquad \text{for all } n\in\mathbb{N}_0,$$
  with $\Lambda_n\in \operatorname{Mat}_N(\mathbb{C})$, then $F_i(x)=\sum_{j=0}^i x^j F_j^i$, $F_j^i \in \operatorname{Mat}_N(\mathbb{C})$, is a polynomial and $\deg(F_i)\leq i$. Moreover, $D$ is determined by the sequence $\{\Lambda_n\}_{n\geq 0}$ and
 \begin{equation}\label{eigenvaluemonicos}
   \Lambda_n=\Lambda_n(D)=\sum_{i=0}^s [n]_i F_i^i, \qquad \text{for all } n\geq 0,
 \end{equation}
    where $[n]_i=n(n-1)\cdots (n-i+1)$, $[n]_0=1$.
\end{prop}

Given a weight matrix $W$, the algebra 
\begin{equation}\label{algDW}
  \mathcal D(W)=\{D\in \operatorname{Mat}_{N}(\Omega[x]) : \, P_nD=\Lambda_n(D) P_n, \, \Lambda_n(D)\in \operatorname{Mat}_N(\mathbb{C}), \text{ for all }n\in\mathbb{N}_0\}
\end{equation} 
is introduced in \cite{GT07}, where $\{P_n\}_{n\in \mathbb{N}_0}$ is any sequence of matrix-valued orthogonal polynomials with respect to $W$.
We observe that the definition of $\mathcal D(W)$ depends only on the weight matrix $W$ and not  on the particular sequence of orthogonal polynomials,
since two sequences $\{P_n\}_{n\in\mathbb{N}_0}$ and $\{Q_n\}_{n\in\mathbb{N}_0}$ of matrix orthogonal polynomials with respect to the weight $W$ are related by
$P_n=M_nQ_n$, with $\{M_n\}_{n\in\mathbb{N}_0}$ invertible matrices (see \cite[Corollary 2.5]{GT07}).

\begin{prop} [\cite{GT07}, Proposition 2.8]\label{prop2.8-GT}
For each $n\in\mathbb{N}_0$, the mapping $D\mapsto \Lambda_n(D) $ is a representation of $\mathcal D(W)$ in $\operatorname{Mat}_N(\mathbb{C})$. Moreover, the sequence of representations $\{\Lambda_n\}_{n\in\mathbb{N}_0}$ separates the elements of $\mathcal D(W)$, i.e.
if $\Lambda_n(D_1)=\Lambda_n(D_2)$ for all $n\geq 0$ then $D_1=D_2$. 
\end{prop}

The {\em formal adjoint} on $\operatorname{Mat}_{N}(\Omega[[x]])$, denoted  by $\mbox{}^*$, is the unique involution extending Hermitian conjugate on $\operatorname{Mat}_{N}(\mathbb{C}[x])$ and sending $\partial I$ to $-\partial I$. 
The {\em formal $W$-adjoint} of $ \mathfrak{D}\in \operatorname{Mat}_{N}(\Omega[x])$, or  the formal adjoint of $\mathfrak D$ with respect to $W(x)$  is the differential operator $\mathfrak{D}^{\dagger} \in \operatorname{Mat}_{N}(\Omega[[x]])$ defined
by
\begin{equation}\label{W-daga}
    \mathfrak{D}^{\dagger}:= W(x)\mathfrak{D}^{\ast}W(x)^{-1},
\end{equation}
where  $\mathfrak{D}^{\ast}$ is the formal adjoint of $\mathfrak D$. 
An operator $\mathfrak{D}\in \operatorname{Mat}_{N}(\Omega[x])$ is called {\em $W$-adjointable} if there exists 
$\widetilde {\mathfrak{D}} \in \operatorname{Mat}_{N}(\Omega[x])$, such that
$$\langle P\mathfrak{D},Q\rangle=\langle P,Q\widetilde{\mathfrak{D}}\rangle,$$ for all $P,Q\in \operatorname{Mat}_N(\mathbb{C}[x])$. Then we say that the operator $\widetilde{\mathfrak D}$ is the $W$-adjoint of $\mathfrak D $.

\begin{remark}
In this work, to avoid pathological weights and their technicalities, we consider only weight matrices that are  ``good enough", i.e., such that for each integer $n \geq 0$ the $n$-th derivative $W^{(n)}(x)$ decreases exponentially at infinity and there exists a scalar polynomial $p_{n}(x)$ such that $W^{(n)}(x)p_{n}(x)$ has finite moments. See \cite[Section 2.2]{CY18}.
\end{remark}
 
With the above remark in mind, it follows that the $W$-adjoint of any $W$-adjointable operator in $\operatorname{Mat}_{N}(\Omega[x])$ coincides with the formal $W$-adjoint. 
For $\mathfrak D= \sum_{j=0}^n \partial ^j F_j \in \operatorname{Mat}_{N}(\Omega[x]) $,  the formal $W$-adjoint of $\mathfrak D$ is given by $\mathfrak{D}^\dagger= \sum_{k=0}^n \partial ^k G_k$, with
\begin{equation}\label{daga}
    G_k=  \sum_{j = 0}^{n-k}(-1)^{n-j} \binom{n-j}{k}
    (WF_{n-j}^*)^{(n-k-j)} W^{-1} , \qquad \text{for } 0\leq k\leq n.  
 \end{equation}
 It is a matter of careful integration by parts to see that $\langle \, P\mathfrak{D}, Q\, \rangle=\langle \, P,\, Q {\mathfrak{D}}^\dagger \rangle$ if  the following set of ``boundary conditions" are satisfied
  \begin{equation} 
  \lim_{x\to x_i}\, \sum_{j = 0}^{p-1}(-1)^{n-j+p-1} \binom{n-j}{k}
    \big (F_{n-j}(x)W(x)\big)^{(p-1-j)}=0, 
      \end{equation}
for $1\leq p\leq n$  and $0\leq k\leq n-p$, where $x_i$ are the endpoints of the support of the weight $W$.

\

We say that  a differential operator $D\in \mathcal D(W)$ is $W$-{\em symmetric} if $\langle P{D},Q\rangle=\langle P,QD\rangle$, for all $P,Q\in \operatorname{Mat}_N(\mathbb{C}[x])$.
In particular, since all operators in $\mathcal{D}(W)$ are $W$-adjointable, we have that $D \in \mathcal{D}(W)$ is $W$-symmetric if and only if $D = D^{\dagger}$.

The condition of symmetry for a differential operator in the algebra $\mathcal D(W)$ is equivalent to the following set of differential equations involving the weight $W$ and the coefficients of  $D$.

\begin{thm}\label{equivDsymm}
 Let $\mathfrak{D} =\sum_{i=0}^n \partial^i F_i(x)$ be a  differential operator of order $n$ in $\mathcal{D}(W)$.
 Then $\mathfrak{D}$ is $W$-symmetric if and only if
\begin{equation*}
    \sum_{j = 0}^{n-k}(-1)^{n-j} \binom{n-j}{k}
    (F_{n-j}W)^{(n-k-j)}= WF_{k}^\ast   
\end{equation*}
for all $0\leq k \leq n$.
 \end{thm}
 \begin{proof}  
 An operator $\mathfrak{D}\in \mathcal{D}(W)$ is $W$-symmetric if and only if $\mathfrak{D}= \mathfrak{D}^\dagger$. The statement follows by using the explicit expression of the coefficients of  $\mathfrak{D}^\dagger$ given in  \eqref{daga}. 
 \end{proof}

In particular, the coefficients of a differential operator of order two in $\mathcal S(W)$, the set of all $W$-symmetric operators in $\mathcal D(W)$, satisfy the classical {\em symmetry equations}: 
\begin{equation}\label{symmeq2}
\begin{split}
  F_2 W &=WF_2^*,  \\
    2(F_2W)'-F_1W &=WF_1^*,\\
    (F_2W)''-(F_1W)'+F_0W & =WF_0^*.   
\end{split}   
\end{equation}

\smallskip

On the other hand, in \cite[Corollary 4.5]{GT07} it is proved that the algebra $\mathcal D(W)$ satisfies

$$\mathcal D(W)= \mathcal S (W)\oplus i \mathcal S (W).$$

\begin{definition}
   We say that a differential operator  
     $\mathcal{V} \in \operatorname{Mat}_{N}(\Omega[x])$  is a  
     \textbf{degree-preserving} operator if the degree of $P(x) \cdot \mathcal{V}$ is equal to the degree of $P(x)$ for all $P(x) \in \operatorname{Mat}_{N}(\mathbb{C}[x])$.
\end{definition}

In particular, a degree-preserving differential operator $\mathcal{V} = \sum_{j=0}^{m}\partial^{j}F_{j} \in \operatorname{Mat}_{N}(\Omega[x])$ satisfies $\deg(F_{j})\leq j$ for all $j$. If $P(x) \in \operatorname{Mat}_{N}(\mathbb{C}[x])$ has a nonsingular leading coefficient, then $P(x)\cdot \mathcal{V}$ also has a nonsingular leading coefficient. Furthermore, it should be noted that every degree-preserving differential operator is not a zero divisor in the algebra of differential operators $\operatorname{Mat}_{N}(\Omega[x])$.

\section{Strong Darboux transformation}
In this section, we study Darboux transformations between weight matrices $W$ and $\widetilde W$ and the relationship between their algebras $\mathcal{D}(W)$ and $\mathcal D(\widetilde W)$.  From \cite{BP23b}, we have the following definition of Darboux transformation.
\begin{definition}\label{darb def}
    Let $W$ and $\widetilde{W}$ be weight matrices supported on the same interval, and let $P_{n}(x)$ and $\widetilde{P}_{n}(x)$ be their associated sequence of monic orthogonal polynomials. We say that $\widetilde{W}$ is a \textbf{Darboux transformation} of $W$ if there exists a differential operator $D \in \mathcal{D}(W)$ that can be factorized as $D = \mathcal{V}\mathcal{N}$ with degree-preserving operators $\mathcal{V}, \, \mathcal{N} \in \operatorname{Mat}_{N}(\Omega[x])$, such that
    \begin{equation*}
            P_{n}(x) \cdot \mathcal{V}  = A_{n}\widetilde{P}_{n}(x)
    \end{equation*}
for a sequence of matrices $A_{n} \in \operatorname{Mat}_{N}(\mathbb{C})$.
\end{definition}

We proved in \cite{BP23b} that the Darboux transformation defines an equivalence relation on $\mathbb{W}(N,\mathcal{I})$, the set of all weight matrices of size $N$ supported on the interval $\mathcal{I}$. 
We also proved the following results. 

\begin{prop} \label{cons}
Let $W$ and $\widetilde{W}$ be weight matrices and let $P_{n}(x)$ and $\widetilde{P}_{n}(x)$ be their associated sequence of orthogonal monic polynomials. If  $\widetilde{W}$ is a Darboux transformation of $W$ with $D = \mathcal{V}\mathcal{N}$ as in Definition \ref{darb def}, then 
\begin{enumerate}[(i)] 
    \item $P_{n}(x)\cdot\mathcal{V}$ is a sequence of orthogonal polynomials for $\widetilde{W}$,
    \item $\widetilde{D} = \mathcal{N}\mathcal{V}$ belongs to $\mathcal{D}(\widetilde{W})$ and $W$ is a Darboux transformation of $\widetilde{W},$
    \item the eigenvalues $\Lambda_{n}(D)$ are nonsingular for all $n\geq 0$.
\end{enumerate}
\end{prop}

When  $W$ is a Darboux transformation of $\widetilde W$ the algebras $\mathcal{D}(W)$ and $\mathcal D(\widetilde W) $ are closely related to each other.   

\begin{prop}\label{alg debil}
     Let $\widetilde{W}, \, W$ be weight matrices of size $N$ such that $\widetilde{W}$ is a Darboux transformation of $W$ with $D = \mathcal{V}\mathcal{N} \in \mathcal{D}(W)$, as in  Definition \ref{darb def}. Then the algebras satisfy
    $$\mathcal{V}\mathcal{D}(\widetilde{W})\mathcal{N} \subseteq \mathcal{D}(W) \quad \text{ and } \quad \mathcal{N} \mathcal{D}(W)\mathcal{V} \subseteq \mathcal{D}(\widetilde{W}).$$
\end{prop}

The following proposition gives us some conditions that guarantee that $\widetilde{W}$ is a Darboux transformation of $W$.

\begin{prop}\label{strong1}
      Let $W$ and  $\widetilde{W}$  be weight matrices of size $N$. If there exist a degree-preserving differential operator $\mathcal{V}$ such that $D = \mathcal{V}\mathcal{N} \in \mathcal{D}(W)$, with $\mathcal{N} = \widetilde{W}\mathcal{V}^{\ast}W^{-1}$, and
    $$\langle P\cdot \mathcal{V},Q \rangle_{\widetilde{W}} = \langle P, Q\cdot \mathcal{N} \rangle_{W} \qquad \text{ for all } P,Q \in \operatorname{Mat}_{N}(\mathbb{C}[x]),$$ then $\widetilde{W}$ is a Darboux transformation of $W$.  
\end{prop}
\begin{proof}
Let $P_{n}(x)$ and $\widetilde{P}_{n}(x)$ be the associated sequence of monic orthogonal polynomials for $W$ and $\widetilde{W}$ respectively.
    We need to show that $P_{n}(x)\cdot \mathcal{V} = A_{n}\widetilde{P}_{n}(x)$ for some sequence of matrices $A_{n} \in \operatorname{Mat}_{N}(\mathbb{C})$, and that $\mathcal{N}$ is a degree-preserving operator. 
    Since $\mathcal{V}$ is degree-preserving, we have that  $P_{n}(x)\cdot \mathcal{V}$ is a polynomial of degree $n$ whose leading coefficient is a nonsingular matrix. If $n\not=m$,  then 
\begin{align*}
        \langle P_{n}(x)\cdot \mathcal{V}, P_{m}(x)\cdot \mathcal{V} \rangle_{\widetilde{W}} & = \langle P_{n}(x), P_{m}(x)\cdot \mathcal{V}\mathcal{N} \rangle_{W}
        = \langle P_{n}(x), P_{m}(x)\cdot D \rangle_{W} \\ & 
        = \langle P_{n}(x), P_{m}(x) \rangle_{W} \Lambda_{m}(D)^{\ast} = 0.
    \end{align*}
Therefore $P_{n}(x)\cdot \mathcal{V}$ is a sequence of orthogonal polynomials for $\widetilde{W}$. Thus, there exists a sequence of nonsingular matrices $A_{n}$ such that $A_{n}\widetilde{P}_{n}(x) = P_{n}(x) \cdot \mathcal{V}$. Finally, to demonstrate that $\mathcal{N}$ is a degree-preserving operator, consider $Q \in \operatorname{Mat}_{N}(\mathbb{C}[x])$ of degree $m$. 
We can express the polynomial $Q$ as $Q(x) = \sum_{j=0}^{m}K_{j}\widetilde{P}_{j}(x)$, where $K_{j} \in \operatorname{Mat}_{N}(\mathbb{C})$ and $K_{m} \neq 0$. Now, if $Q(x)\cdot \mathcal{N}$ is of degree lower than $m$, then 
$$A_{m}\|\widetilde{P}_{m}\|^{2}K_{m} = \langle P_{m}\cdot\mathcal{V}, Q \rangle_{\widetilde{W}} = \langle P_{m}, Q\cdot \mathcal{N} \rangle_{W} = 0,$$
which is a contradiction.
\end{proof}

The above specific factorization of a differential operator $D \in \mathcal{D}(W)$ as $ D=\mathcal V \widetilde{W}\mathcal{V}^{\ast}W^{-1} $ is very powerful. We can obtain the sequence of orthogonal polynomials for $\widetilde{W}$ knowing that of $W$. Moreover, we will see that with this factorization we can map symmetric operators to symmetric ones. This motivates us to introduce the following definition.

\begin{definition}\label{strong darb}
    Let $W$, $\widetilde{W} \in \mathbb{W}(N,\mathcal{I})$. We say that $\widetilde{W}$ is a \textbf{strong Darboux transformation} of $W$ if there exist a degree-preserving differential operator $\mathcal{V}$ with nonsingular leading coefficient such that $D = \mathcal{V}\mathcal{N} \in \mathcal{D}(W)$ with $\mathcal{N} = \widetilde{W}\mathcal{V}^{\ast}W^{-1}$ a differential operator satisfying $$\langle P\cdot \mathcal{V},Q \rangle_{\widetilde{W}} = \langle P, Q\cdot \mathcal{N} \rangle_{W}, \qquad \text{  for all } P,Q \in \operatorname{Mat}_{N}(\mathbb{C}[x]).$$ We call $\mathcal{V}$ a \textbf{strong Darboux transformer} from $W$ to $\widetilde{W}$.
\end{definition}

 \begin{remark}
If $\widetilde W$ is a  strong Darboux transformation of $W$, then it is also a Darboux transformation, (see  Proposition \ref{strong1}). However, there are weights $\widetilde W$ and $W$ that are Darboux transformations but they are not strong Darboux transformations, as we will show in Examples \ref{example} and \ref{example 2}.
 \end{remark}

\smallskip
Given two weight matrices $W$ and $\widetilde{W}$  supported on the same interval $\mathcal I$ we introduce the relation 
\begin{equation*}\label{equivrelat}
W \sim \widetilde{W}  \qquad \text{ if and only if }  \qquad \widetilde{W}  \text{ is a strong  Darboux transformation of $W$}. 
\end{equation*}

\begin{prop}
 The strong Darboux transformation defines an equivalence relation on $\mathbb{W}(N,\mathcal{I})$. 
\end{prop}
\begin{proof}
The reflexivity follows by taking the strong Darboux transformer $\mathcal{V} = I$.

Suppose that $\widetilde{W}$ is a strong Darboux transformation of $W$, with strong Darboux transformer $\mathcal{V}$ and $\mathcal{N} = \widetilde{W}\mathcal{V}^{\ast}W^{-1}$ as in the definition. The leading coefficient of $\mathcal{N}$ is nonsingular because the one of $\mathcal{V}$ is. By Proposition \ref{cons} it follows that $\widetilde{D} = \mathcal{N}\mathcal{V} \in \mathcal{D}(\widetilde{W})$. Besides, 
$W\mathcal{N}^{\ast}\widetilde{W}^{-1}= \mathcal{V}$ and  $\langle P\cdot \mathcal{N}, Q \rangle_{W} = (\langle Q, P\cdot \mathcal{N}\rangle_{W})^{\ast} = (\langle Q\cdot \mathcal{V}, P \rangle_{\widetilde{W}})^{\ast} = \langle P, Q\cdot \mathcal{V} \rangle_{\widetilde{W}}$ for all $P,Q\in \operatorname{Mat}_{N}(\mathbb{C}[x])$. Therefore  $W$ is a strong Darboux transformation of $\widetilde{W}$ with strong Darboux transformer $\mathcal{N}$.

 Finally, let $W_{1}$, $W_{2}$ and $W_{3} 
 $ be weight matrices  such that $W_1 \sim W_{2}$ and $W_2\sim W_3$. Let $\mathcal{V}_{1}$ be a strong Darboux transformer from $W_{1}$ to $W_{2}$ and  $\mathcal{V}_{2}$ a strong Darboux transformer from $W_{2}$ to $W_{3}$. It is easy to verify that the differential operator $\mathcal{V}_{1}\mathcal{V}_{2}$ is a strong Darboux transformer from $W_{1}$ to $W_{3}$ and therefore  $W_1\sim W_{3}$ and this completes the proof.     
\end{proof}

\begin{remark}
    If   $W,\widetilde{W} \in \mathbb{W}(N,\mathcal I)
    $ are equivalent weights, that is if there exists an invertible matrix $M \in \operatorname{Mat}_{N}(\mathbb{C})$ such that $\widetilde{W} = MWM^{\ast}$ then $\widetilde{W}$ is a strong Darboux transformation of $W$ with strong Darboux transformer $\mathcal{V} = M^{-1}$. 
\end{remark}

The following theorem is the main tool for acquiring information about the algebra $\mathcal D(W)$ by understanding the algebra $\mathcal D(\widetilde W)$ when $W$ and $\widetilde W$ are strong Darboux transformations of each other.

Recall that $\mathcal S(W)$ denotes the set of all $W$-symmetric differential operators in $\mathcal D(W)$ and we have $\mathcal D(W)= \mathcal S (W)\oplus i \mathcal S (W)$ as real vector spaces.  

\begin{thm}\label{main}
Let $\widetilde{W}$ be a strong Darboux transformation of $W$. 
If  $\mathcal{V}$ is a strong Darboux transformer from $W$ to $\widetilde{W}$ with $\mathcal{N} =  \widetilde{W}\mathcal{V}^{\ast}W^{-1}$ and $D =\mathcal{V}\mathcal{N}\in \mathcal D(W)$, then  $$\mathcal{V}\mathcal{S}(\widetilde{W})\mathcal{N} \subset \mathcal{S}(W). $$ 
\end{thm}
\begin{proof} 
     If  $\mathcal{A} \in \mathcal{S}(\widetilde{W})$ then, by Proposition \ref{alg debil}, we have that $\mathcal{B} = \mathcal{V}\mathcal{A}\mathcal{N}$ is in $\mathcal{D}(W)$. 
     Now 
     $$\mathcal{B}^{\dagger} = W(\mathcal{N}^{\ast}\mathcal{A}^{\ast}\mathcal{V}^{\ast})W^{-1} = W(W^{-1}\mathcal{V}\widetilde{W}\mathcal{A}^{\ast}\mathcal{V}^{\ast})W^{-1}=\mathcal{V}\mathcal{A}^{\dagger}\mathcal{N} = \mathcal{V}\mathcal{A}\mathcal{N}=\mathcal{B}.  $$
     Thus $\mathcal B$ is a symmetric operator in $\mathcal D(W)$.
\end{proof}

In the next theorem, we see that the strong Darboux transformation preserves the parity of the operators. In particular, if all differential operators in $\mathcal{D}(W)$ have even order as in the case of $W(x) = w(x)I$ for some classical scalar weight $w$, then all differential operator in $\mathcal{D}(\widetilde{W})$ has even order.

\begin{thm} \label{parity}
Let $\widetilde{W}$ be a strong Darboux transformation of $W$ and let $\mathcal{V}$ be a strong Darboux transformer from $W$ to $\widetilde{W}$ with $\mathcal{N} = \widetilde{W}\mathcal{V}^{\ast}W^{-1}$.
If $\widetilde{D} \in \mathcal{D}(\widetilde{W})$ is of even (odd) order, then $\mathcal{D} = \mathcal{V}\widetilde{D}\mathcal{N}\in \mathcal D(W)$ is of even (odd) order.
\end{thm}

\begin{proof}
Let $\widetilde{D} \in \mathcal{D}(\widetilde{W})$ and $\mathcal{D} = \mathcal{V}\widetilde{D}\mathcal{N}$. Since $\mathcal{V}$ has nonsingular leading coefficient, then $\mathcal{N} = \widetilde{W}\mathcal{V}^{\ast}W^{-1}$ has also nonsingular leading coefficient and $\ell = \operatorname{Ord}(\mathcal{V}) = \operatorname{Ord}(\mathcal{N})$. Therefore it follows that 
$ \operatorname{Ord}(\mathcal{D}) = \operatorname{Ord}(\mathcal{V}\widetilde{D}\mathcal{N}) = \operatorname{Ord}(\mathcal{V}) + \operatorname{Ord}(\widetilde{D}) + \operatorname{Ord}(\mathcal{N}) =  2\ell + \operatorname{Ord}(\widetilde{D})$.
Thus $\operatorname{Ord}(\mathcal{D})$ has the same parity of  $\operatorname{Ord}(\widetilde{D})$. 
\end{proof}

\begin{cor}\label{parity cor}
    Let $\widetilde{W}$ be a strong Darboux transformation of $W$. If every differential operator in $\mathcal{D}(W)$ is of even order, then every differential operator in $\mathcal{D}(\widetilde{W})$ will also be of even order.
\end{cor}
\begin{proof}
    The statement follows from Theorem \ref{parity}  and Proposition \ref{alg debil}.
\end{proof}

Now, we present a straightforward example that illustrates how a weight matrix $\widetilde W$ can be a Darboux transformation of a weight $W$, yet fail to be a strong Darboux transformation.

\begin{example}\label{example}
    Let us consider the Laguerre weight $w_{\alpha}(x) = e^{-x}x^{\alpha}$, for  $\alpha > -1$. 
    
    The matrix weight 
 $w_{\alpha+1}(x)\oplus w_{\alpha}(x)$ is a Darboux transformation of $w_{\alpha}(x) \oplus w_{\alpha}(x)$ but it is not  a strong Darboux transformation of it. 
    
   In fact, we have proved in \cite{BP23b} that $w_{\alpha+1}$ is a Darboux transformation of $w_{\alpha}$. Then it follows that $w_{\alpha+1}(x) \oplus w_{\alpha}(x)$ is a Darboux transformation of $ w_{\alpha}(x) \oplus w_{\alpha}(x)$.  (See Proposition 3.7 in \cite{BP23b}).   
    Also from \cite{BP23b} we have that 
$$\mathcal D(w_\alpha\oplus w_\alpha )= \begin{pmatrix} \mathcal{D}(w_{\alpha}) && \mathcal{D}(w_{\alpha}) \\ \mathcal{D}(w_{\alpha}) && \mathcal{D}(w_{\alpha}) \end{pmatrix}= \begin{pmatrix}
    \mathbb C[\delta] & \mathbb C[\delta] \\ \mathbb C[\delta] & \mathbb C[\delta] 
\end{pmatrix},$$
 where  $ \delta=\partial^2 x + \partial(\alpha+1-x)$. In particular, all differential operators in $\mathcal{D}(w_{\alpha}\oplus w_{\alpha})$ are of even order. 
 On the other hand, we have that the differential operator $\begin{pmatrix} 0 && 0 \\ \partial - 1 && 0 \end{pmatrix}$ belongs to the algebra $\mathcal{D}(w_{\alpha+1}\oplus w_{\alpha})$, (see Proposition 5.4 in \cite{BP23b}). 
    Therefore, by Corollary \ref{parity cor}, we have that 
    $w_{\alpha+1}(x)\oplus w_{\alpha}(x)$ is not  a strong Darboux transformation of $w_{\alpha}(x) \oplus w_{\alpha}(x)$. 
\end{example}

Now we present a much more interesting example. We construct an irreducible $3 \times 3$ weight matrix of Hermite-type, which is a Darboux transformation of $e^{-x^{2}}I$, but it is not a strong Darboux transformation of it.

\begin{example} \label{example 2}
    Let $W$ be the weight matrix 
    \begin{equation}\label{wz}
        W(x) = e^{-x^{2}}\begin{pmatrix} a^{2}x^{2} + 1 && ax && 0 \\ ax && b^{2}x^{2} + 1 && bx \\ 0 && bx && 1 \end{pmatrix}, \quad a,b \in \mathbb{R}-\{0\}, \, x \in \mathbb{R}.
    \end{equation}
    \begin{prop}\label{pz}
        $W$ is a Darboux transformation of $\widetilde{W}(x) = e^{-x^{2}}I$ but not a strong Darboux transformation of it.
    \end{prop}
    \begin{proof}
        We have proved in \cite{BP23b} that the algebra $\mathcal{D}(\widetilde{W})$, for $\widetilde{W}(x)=e^{-x^{2}}I$, is 
        $$\mathcal D(\widetilde{W})=\sum_{j=1}^{3}\sum_{i=1}^{3} \mathbb{C}[\delta]E_{i,j},$$ where $\delta = \partial^{2} + \partial(-2x)$ is the second-order differential operator of the classical scalar Hermite weight, and $E_{i,j}$ denotes the matrix with zeros everywhere except for a $1$ in the $(i,j)$-entry. Then, the operator
        $$D = \begin{pmatrix} - \delta + 2 \frac{(a^{2}+2)}{a^{2}} && - \frac{2}{a}\delta + 4 \frac{(a^{2}+2)}{a^{3}} && 0 \\ -\frac{2}{a}\delta + 4 \frac{(a^{2}+2)}{a^{3}} && \delta^{2} + 8 \frac{(a^{2}b^{2}+2(a^{2} + b^{2}))}{a^{4}b^{2}} && 0 \\ 0 && 0 && - \delta^{3} + 2 \frac{(a^{2}+4)}{a^{2}}\delta^{2} - 8 \frac{(a^{2}+2)}{a^{4}}\delta \end{pmatrix}\in \mathcal{D}(e^{-x^{2}}I) .$$
        
        We can factorize $D$ as $D = \mathcal{V}\mathcal{N}$ with 
        \begin{align*}
                  \mathcal{V} & = \partial^{3} \begin{psmallmatrix} 0 && 0 && 0 \\ 0 && 0 && 0 \\ 0 && 0 && 1 \end{psmallmatrix} + \partial^{2} \begin{psmallmatrix} 0 && 0 && 0 \\ 1 && -ax && abx^{2} \\ 0 && 0 && -4x\end{psmallmatrix} + \partial \begin{psmallmatrix} 0 && 1 && -bx \\ 0 && 0 && 0 \\ -\frac{4}{ab} && \frac{4x}{b} && -2 \frac{(a^{2}+2)}{a^{2}} \end{psmallmatrix} + \begin{psmallmatrix} -\frac{2}{a} && 0 && 0 \\- \frac{4}{a^{2}} && 0 && \frac{4}{ab} \\ 0 && \frac{8}{ba^{2}} && 0 \end{psmallmatrix}, \\ \intertext{and }
                  \mathcal{N} & = \partial^{3} \begin{psmallmatrix}0 && 0 && 0 \\ 0 && 0 && -bx \\ 0 && 0 && -1 \end{psmallmatrix} + \partial^{2} \begin{psmallmatrix} 0 && 1 && 0 \\ 0 && 0 && b(2x^{2} - 3) \\0 && 0 && 2x\end{psmallmatrix} + \partial \begin{psmallmatrix} -ax && -4x && \frac{4}{ab} \\ - 1 && 0 && 4bx\frac{(a^{2}+1)}{a^{2}} \\ 0 && 0 && \frac{4}{a^{2}}\end{psmallmatrix} + \begin{psmallmatrix} - \frac{(a^{2}+2)}{a} && -2 \frac{(a^{2}+2)}{a^{2}} && 0 \\ 0 && 0 && 4 \frac{(b^{2}+2)}{ba^{2}} \\ 0 && \frac{4}{ab} && 0\end{psmallmatrix}.
        \end{align*}
        The differential operators $\mathcal{V}$ and $\mathcal{N}$ are degree-preserving, and it follows that $\mathcal{N} = W\mathcal{V}^{\ast}\widetilde{W}^{-1}$. The equality $\langle P \cdot \mathcal{V}, Q \rangle_{W} = \langle P, Q\cdot \mathcal{N} \rangle_{\widetilde{W}}$ holds for all $P, \, Q \in \operatorname{Mat}_{3}(\mathbb{C}[x])$ by integrating by parts since the boundary conditions will vanish due to the exponential decay of $e^{-x^{2}}$. Thus, by Proposition \ref{strong1} we have that $W$ is a Darboux transformation of $e^{-x^{2}}I$.

        Now, we assume that $W$ is a strong Darboux transformation $\widetilde{W}(x) = e^{-x^{2}}I$. Since $\mathcal{D}(e^{-x^{2}}I)$ contains only even-order differential operators, we have, by Corollary \ref{parity cor}, that $\mathcal{D}(W)$ must contain only even-order differential operators. However, by Proposition \ref{alg debil}, we have that $\mathcal{A} = \mathcal{N} \begin{psmallmatrix} 0 && 0 && 0 \\ 0 && 1 && 0 \\ 0 && 0 && 0 \end{psmallmatrix} \mathcal{V} \in \mathcal{D}(W)$ and is of odd order, which is a contradiction.
    \end{proof}
\end{example}

\begin{remark}
    As a consequence of Proposition \ref{pz}, we obtain an explicit expression of a sequence of orthogonal polynomials associated to the weight $W$ given in \eqref{wz}. The sequence is given by 
    $$Q_{n}(x) =H_{n}(x)\cdot \mathcal{V} = \begin{psmallmatrix} -\frac{2}{a} H_{n}(x) & H_{n}'(x) & -bxH_{n}'(x) \\ H_{n}''(x)-\frac{4}{ab}H_{n}(x) & -axH_{n}''(x) & abx^{2}H_{n}''(x)+\frac{4}{ab}H_{n}(x) \\ -\frac{4}{ab} H_{n}'(x) & \frac{4x}{b} H_{n}'(x) + \frac{8}{ba^{2}} H_{n}(x) & H_{n}(x)'''-4xH_{n}''(x)-\frac{2(a^{2}+2)}{a^{2}} H_{n}'(x)\end{psmallmatrix},$$
    where $H_{n}(x)$ is the sequence of monic orthogonal polynomials for the classical scalar weight of Hermite.
\end{remark}

\section{The algebra $\mathcal{D}(W)$ for a Laguerre example}\label{Laguerre-Sect}

In this section, we illustrate the efficacy of the strong Darboux transformations by computing the algebra $\mathcal{D}(W)$. We consider the irreducible 
Laguerre-type weight matrix 
\begin{equation}\label{Lag weight}
    W(x) = e^{-x}x^{\alpha}\begin{pmatrix} x(1+a^{2}x) && ax \\ax && 1 \end{pmatrix}, \quad \alpha > -1, \, a\not=0, \, x\in (0,\infty).
\end{equation}
 This weight was introduced and studied in \cite{DI08}, where based on computational evidence, the authors conjectured the generators of the algebra $\mathcal D(W)$.

Here, we establish that the weight $W$ is a strong Darboux transformation of a direct sum of classical scalar weights. Consequently, by applying the theorems within this framework, we determine the algebra $\mathcal D(W)$ and identify its generators, thus proving the original conjecture. Additionally, we provide explicit expressions for a sequence of orthogonal polynomials in terms of classical scalar Laguerre weights.
 
    \begin{prop} \label{ddd}
    The irreducible weight matrix $W(x) = e^{-x}x^{\alpha}\begin{pmatrix} x(1+a^{2}x) && ax \\ax && 1 \end{pmatrix}$ is a strong Darboux transformation of the direct sum of scalar Laguerre weights $\widetilde W(x) = \begin{pmatrix} e^{-x}x^{\alpha+1} && 0 \\ 0 && e^{-x}x^{\alpha} \end{pmatrix}$.
    \end{prop}
    \begin{proof}
    The differential operator 
    $$D = -\partial^{2} xI - \partial \begin{pmatrix} \alpha + 2 -x && 0 \\ 0 && \alpha+1-x \end{pmatrix} + \begin{pmatrix}1+\frac{1}{a^{2}} && 0 \\ 0 && \frac{1}{a^{2}} \end{pmatrix} \in \mathcal{D}(\widetilde W)$$
    can be  factorized as $D = \mathcal{V}\mathcal{N}$, with the degree-preserving differential operators
    $$\mathcal{V} = \partial \begin{pmatrix} 0 && x \\ -1 && ax \end{pmatrix} + \begin{pmatrix} -\frac{1}{a} && \alpha + 1 \\ 0 && \frac{1}{a} \end{pmatrix},$$ and 
    $$\mathcal{N} = W\mathcal{V}^{\ast}\widetilde{W}^{-1} = \partial \begin{pmatrix} -ax && x \\-1 && 0 \end{pmatrix} + \begin{pmatrix} -a-\frac{1}{a} && \alpha+1 \\ 0 && \frac{1}{a} \end{pmatrix}.$$
     Now, we need to prove that  $\langle P\cdot \mathcal{V}, Q \rangle_{W} = \langle P , Q\cdot \mathcal{N} \rangle_{\widetilde W}$, for all $P, \, Q\in\operatorname{Mat}_{2}(\mathbb{C}[x])$. We observe that $W(x) = T(x)\widetilde W(x)T(x)^{\ast}$ with $T(x) = \begin{pmatrix}1 && ax \\ 0 && 1 \end{pmatrix}$ and 
    \begin{equation*}
        \begin{split}
            \langle P\cdot \mathcal{V}, Q \rangle_{W} & 
            = \int_{0}^{\infty}P(x)\cdot(\mathcal{V}T(x))\widetilde W(x)(Q(x)T(x))^{\ast}dx  
            = \langle P\cdot (\mathcal{V}T), QT\rangle_{\widetilde W}.
        \end{split}
    \end{equation*}
It is a matter of integration by parts, to obtain the differential operator 
$\mathcal{V}T = \partial \begin{pmatrix} 0 && x \\ -1 && 0 \end{pmatrix} + \begin{pmatrix} -\frac{1}{a} && -x + \alpha + 1 \\ 0 && \frac{1}{a} \end{pmatrix}$ is $\widetilde W$-adjointable because $\lim_{x\to 0, \infty} P(x)\begin{pmatrix} 0 && x \\ -1 && 0 \end{pmatrix} \widetilde W(x) Q(x)^{\ast} = 0$ for all $P, \, Q \in \operatorname{Mat}_{2}(\mathbb{C}[x])$.
    The  $\widetilde W$-adjoint of $\mathcal VT(x)$ is  given by
    $$(\mathcal{V}T)^{\dagger} = \widetilde W(\mathcal{V}T)^{\ast}\widetilde W^{-1} = \partial \begin{pmatrix} 0 && x \\ -1 && 0 \end{pmatrix} + \begin{pmatrix} -\frac{1}{a} && \alpha +1 -x \\ 0 && \frac{1}{a} \end{pmatrix}.$$
    We observe that $T (\mathcal{V}T)^{\dagger}=\mathcal N$ and   thus 
    $$\langle P\cdot \mathcal{V}, Q \rangle_{W} = \langle P\cdot (\mathcal{V}T),QT \rangle_{\widetilde W} = \langle P, QT\cdot (\mathcal{V}T)^{\dagger} \rangle_{\widetilde{W}} = \langle P , Q\cdot \mathcal{N} \rangle_{\widetilde W}.$$
    Therefore, $W$ is a strong Darboux transformation of $\widetilde W$. 
    \end{proof}
   
\medskip

As a consequence of Proposition \ref{ddd}, we obtain an explicit expression of the sequence of orthogonal polynomials for $W$ in terms of the sequence of the monic scalar Laguerre orthogonal polynomials. 
\begin{cor}
If   $W$ is the irreducible weight given in \eqref{Lag weight}, then 
$$Q_{n}(x) = \begin{pmatrix} -\frac{1}{a}L_{n}^{\alpha+1}(x) && (\alpha+1)L_{n}^{\alpha+1}(x) + x(L_{n}^{\alpha+1})'(x) \\ -(L_{n}^{\alpha})'(x) && \frac{1}{a}L_{n}^{\alpha}(x) + ax(L_{n}^{\alpha})'(x)\end{pmatrix} , \qquad n\geq 0$$
    is a sequence of orthogonal polynomials for $W$, where $L_{n}^{\alpha}$ is the sequence of monic orthogonal polynomials for the scalar Laguerre weight $w_\alpha(x)=e^{-x}x^{\alpha}$.     
\end{cor}  
\begin{proof}
 The sequence of monic orthogonal polynomials with respect to $\widetilde W(x) = \begin{psmallmatrix} e^{-x}x^{\alpha+1} && 0 \\ 0 && e^{-x}x^{\alpha} \end{psmallmatrix}$ is given by $P_{n}(x) = \operatorname{diag}(L_{n}^{\alpha+1}(x), L_{n}^{\alpha}(x))$.
Hence, from Proposition \ref{cons},  
 $Q_n(x)=P_{n}(x)\cdot \mathcal{V} $ is a sequence of orthogonal polynomials for $W$.
\end{proof}

\medskip 
In order to determine the algebra $\mathcal{D}(W)$, we begin by studying some properties of the algebra $\mathcal{D}(\widetilde W)$, for $\widetilde W(x) = e^{-x}x^{\alpha+1} \oplus e^{-x}x^{\alpha}$. 
We know (see  \cite{BP23b}) that 
    \begin{equation} \label{lag lag}
    \mathcal{D}(\widetilde W) =  \left \{ \begin{pmatrix} p_{1}(\delta_{\alpha+1}) && \mathcal{N}p_{2}(\delta_{\alpha}) \\ \mathcal{V}p_{3}(\delta_{\alpha+1}) && p_{4}(\delta_{\alpha}) \end{pmatrix} \biggm| p_{1}, \, p_{2}, \, p_{3}, \, p_{4} \in \mathbb{C}[x] \right \},
    \end{equation}
    where $\delta_{\alpha} = \partial^{2}x + \partial(\alpha+1-x)$, $\delta_{\alpha+1} = \partial^{2}x + \partial(\alpha+2-x)$, $\mathcal{V} = \partial -1$ and $\mathcal{N} = \partial x + \alpha + 1$. It is noteworthy that all differential operators on the antidiagonal are of odd order, while those on the diagonal are of even order.
    
    The $\widetilde W$-symmetric differential operators in $\mathcal{D}(\widetilde W)$ satisfy the following property. 

    \begin{prop} \label{sim lag}
        Let $D = \sum_{j=0}^{n}\partial^{j}F_{j} \in \mathcal{S}(\widetilde W)$ be a $\widetilde W$-symmetric differential operator or order $n$. If $n$ is even, $n =2m$ for some $m\geq 0$, then the leading coefficient of $D$ satisfies
    $$F_{2m}(x) = \begin{pmatrix} k_{1}x^{m} && 0 \\ 0 && k_{4}x^{m} \end{pmatrix} \text{ for some } k_{1},\,k_{4} \in \mathbb{R}.$$
    If $n$ is odd, $n = 2m + 1$ for some $m\geq 0$, then the leading coefficient of $D$ satisfies
    $$F_{2m+1}(x) = \begin{pmatrix} 0 && (k_{2} + i k_{3})x^{m+1} \\ (-k_{2} + ik_{3})x^{m} && 0 \end{pmatrix} \text{ for some } k_{2}, \, k_{3} \in \mathbb{R}.$$
    \end{prop}
\begin{proof}
    Let  $D$ be a  $\widetilde W$-symmetric operator of order $n$. Since $D=D^\dagger$, from \eqref{daga} with $k=n$,  we get that  
    the leading coefficient of $D$ satisfies $F_{n}(x)\widetilde W(x) = (-1)^{n}\widetilde W(x)F_{n}(x)^{\ast}$. Now the statement follows easily from \eqref{lag lag}.
\end{proof}

\begin{thm} 
The algebra $\mathcal{D}(W)$ 
is generated by $\{ I, D_{1}, D_{2}, D_{3} \}$, where 
    \begin{equation*}
    \begin{split}
        D_{1} & =  \partial^{2} x I + \partial \begin{pmatrix} \alpha + 2 - x && ax \\ 0 && \alpha + 1 -x \end{pmatrix} + \begin{psmallmatrix} -\frac{a^{2}+1}{a^{2}} && a(\alpha+1) \\ 0 && -\frac{1}{a^{2}} \end{psmallmatrix} \\
        D_{2} & = \partial^{2} \begin{psmallmatrix} 0 && -ax^{2} \\ 0 && -x \end{psmallmatrix} + \partial \begin{psmallmatrix} x && -x\frac{(a^{2}\alpha+3a^{2}+1)}{a} \\ \frac{1}{a} &&-\alpha-1 \end{psmallmatrix} + \begin{psmallmatrix} \frac{a^{2}+1}{a^{2}} &&  -\frac{(a^{2}+1)(\alpha+1)}{a} \\  0 &&  0 \end{psmallmatrix}, \\
        D_{3} &= \partial^{3} \begin{psmallmatrix}  ax^{2} && -x^{2}(a^{2}x+1) \\ x &&  -ax^{2} \end{psmallmatrix} + \partial^{2} \begin{psmallmatrix}  ax(\alpha+5)+\frac{2x}{a} && -a^{2}x^{2}(\alpha+5)-x(2\alpha+x+4) \\ 2 + \alpha && -ax(\alpha+2)-\frac{2x}{a} \end{psmallmatrix} \\
         & \quad + \partial  \begin{psmallmatrix} \frac{2(a^{2}+1)(2+\alpha)-x}{a} && -\frac{x}{a^{2}} - 2(a^{2}(\alpha+2)+1)x-(\alpha+2)(\alpha+1) \\ \frac{1}{a^{2}} && \frac{x-2(\alpha+1)}{a} \end{psmallmatrix} + \begin{psmallmatrix} -\frac{\alpha+1}{a} && \frac{(\alpha+1)(a^{2}\alpha-1)}{a^{2}} \\ -\frac{1}{a^{2}} && \frac{\alpha+1}{a} \end{psmallmatrix}.
    \end{split}
\end{equation*}
    
\end{thm}
\begin{proof}
    The weight $W(x)$ is a strong Darboux transformation of $\widetilde W(x) = e^{-x}x^{\alpha+1}\oplus e^{-x}x^{\alpha}$ with strong Darboux transformer $$\mathcal{V} = \partial \begin{pmatrix} 0 && x \\ -1 && ax \end{pmatrix} + \begin{pmatrix} -\frac{1}{a} && \alpha + 1 \\ 0 && \frac{1}{a} \end{pmatrix} \quad \text{ and } \quad  \mathcal{N} =  \partial \begin{pmatrix} -ax && x \\-1 && 0 \end{pmatrix} + \begin{pmatrix} -a-\frac{1}{a} && \alpha+1 \\ 0 && \frac{1}{a} \end{pmatrix}.$$ 
Let us consider the  $\widetilde W$-symmetric differential operators
$$E_{1} = \begin{psmallmatrix} -1 && 0 \\ 0 && -1 \end{psmallmatrix}, \quad  
E_{2} = \begin{psmallmatrix} 1 && 0 \\ 0 && 0 \end{psmallmatrix},  \quad 
E_{3} = \begin{psmallmatrix} 0 && \partial x + \alpha + 1 \\ -\partial +1 && 0 \end{psmallmatrix}, \quad 
E_{4} = \begin{psmallmatrix} 0 && 0 \\ 0 && 1 \end{psmallmatrix}, \quad 
E_{5} = i\begin{psmallmatrix} 0 && \partial x + \alpha + 1 \\ \partial -1 && 0 \end{psmallmatrix}.$$
It is easy to see that $D_{1} = \mathcal{N}E_{1}\mathcal{V}$, $D_{2} = \mathcal{N}E_{2}\mathcal{V}$, and $D_{3} = \mathcal{N}E_{3}\mathcal{V}$. We also define
     \begin{align*}
        D_{4} & = \mathcal{N}E_{4}\mathcal{V} =  \partial^{2} \begin{pmatrix} -x && ax^{2} \\ 0 && 0 \end{pmatrix} + \partial \begin{psmallmatrix}  -\alpha-2 && x\frac{(a^{2}\alpha+2a^{2}+1)}{a} \\  -\frac{1}{a} &&  x \end{psmallmatrix} + \begin{psmallmatrix} 0 && \frac{\alpha+1}{a} \\ 0 &&  \frac{1}{a^{2}}\end{psmallmatrix}, \\ \intertext{and}
         D_{5} & = \mathcal{N}E_{5}\mathcal{V} = \partial^{3} i \begin{psmallmatrix} ax^{2} && -x^{2}(a^{2}x-1) \\  x && -ax^{2}\end{psmallmatrix} + \partial^{2} i \begin{psmallmatrix} a(\alpha+5)x && -x(a^{2}\alpha x+5 a^{2}x-2\alpha+3 x-4) \\ \alpha+2 &&  -a(\alpha+2)x) \end{psmallmatrix} \\ 
         & \quad + \partial i\begin{psmallmatrix} \frac{1}{a}(2a^{2} \alpha+4a^{2}+x) && \frac{-x}{a^{2}} - (\alpha+2)(2a^{2}x-\alpha+4x-1) \\  -\frac{1}{a^{2}} && -\frac{x}{a}\end{psmallmatrix} + i\begin{psmallmatrix} \frac{\alpha+1}{a} && -\frac{1}{a^{2}}(\alpha+1)(a^{2} \alpha+2a^{2}+1) \\  \frac{1}{a^{2}} &&  -\frac{\alpha+1}{a} \end{psmallmatrix}.
    \end{align*}

\noindent From Theorem \ref{main},  we obtain that  $D_{1}$, $D_{2}$,  $D_{3}$, $D_4$ and $D_5$ are $W$-symmetric operators in $\mathcal{D}(W)$.
 
 \smallskip
 Let $E=\sum_{j=0}^{n}\partial^{j}F_{j}(x)$ be a differential operator in $\mathcal S(W)$ of order $n$.  
By Theorem \ref{main} we have that $\mathcal{V}E\mathcal{N} \in \mathcal{S}(\widetilde W)$.
 
 If $n$ is even, $n = 2m$ for some $m\geq 0$, then $\mathcal{V}E\mathcal{N}$ is of order $2m +2$ and by Proposition \ref{sim lag} its leading coefficient is $\begin{pmatrix} k_{1}x^{m+1} && 0 \\ 0 && k_{4}x^{m+1} \end{pmatrix}$ for some $k_{1}, \, k_{4} \in \mathbb{R}$. Hence  the leading coefficient of $E$ is
 \begin{equation} \label{even}
      F_{2m}(x) = k_{4}x^{m-1} \begin{pmatrix} -x && ax^{2} \\ 0 && 0 \end{pmatrix} + k_{1}x^{m-1} \begin{pmatrix} 0 && -ax^{2} \\ 0 && -x \end{pmatrix}. 
\end{equation}
 If $n$ is odd, $n = 2m + 1$ for some $m \geq 0$, then $\mathcal{V}E\mathcal{N}$ is of order $2m + 3$ and by Proposition \ref{sim lag} its leading coefficient is $\begin{pmatrix} 0 && (k_{2} + i k_{3})x^{m+2} \\ (-k_{2} + ik_{3})x^{m+1} && 0 \end{pmatrix}$ for some $k_{2}, \, k_{3} \in \mathbb{R}$. Therefore the leading coefficient of $E$ is 
 \begin{equation}\label{odd}
    F_{2m+1}(x) = k_{2}x^{m-1}\begin{pmatrix}ax^{2} && -x^{2}(a^{2}x+ 1)\\ x && -ax^{2}  \end{pmatrix} + k_{3}x^{m-1}i\begin{pmatrix}ax^{2} && -x^{2}(a^{2}x-1) \\ x && -ax^{2} \end{pmatrix}.
\end{equation}

\ 

   Now we will prove that  the algebra $\mathcal{D}(W)$ is generated by $\mathcal M=\{ I,D_{1}.D_{2},D_{3},D_{4},D_{5} \}$. 
    Let $E\in \mathcal S(W)$ be a differential operator of order $n$. For $n = 0$, we get that $E = k_{1}I$ for some $k_{1} \in \mathbb{R}$. 
    
    There are no first-order differential operators in $\mathcal{D}(W)$, as none satisfy the symmetry equations given in \eqref{symmeq2}. We proceed by induction on $n$ and consider $n>1$.

    For $n=2m$, using \eqref{even}, we have that  the leading coefficient of $E$ is a linear combination of the leading coefficients of $D_{2}$ and $D_{4}$ multiplied by $x^{m-1}$.  Thus the differential operator 
   $\mathcal{A} = E - D_{1}^{m-1}\big (k_{1}D_{2} + k_{4}D_{4}\big )$ has order less than $n$, hence 
    by the inductive hypothesis, the operator $E$ is generated by $\mathcal{M}$. 
    For  $n = 2m+1$, the leading coefficient of $E$ is a linear combination of the leading coefficients of $D_{3}$ and $D_{5}$ multiplied by $x^{m+1}$, see  \eqref{odd}.  Thus the differential operator 
    $ E - D_{1}^{m-1} \big(k_{2}D_{3} + k_{3}D_{5} \big)$ is of order less than $2m+1$ and by inductive hypothesis it follows that $E$ is generated by $\mathcal{M}$.
    Thus, $\mathcal{D}(W)$ is generated by $\{I, D_{1}, D_{2}, D_{3}, D_4,D_5 \}$.

    Finally, we observe that  $D_{4} = - D_{1} - D_{2}$ and by using Proposition \ref{prop2.8-GT} and the expressions of the eigenvalues $\Lambda_n(D)$ 
    given in Proposition \ref{eigenvalue-prop}, it is easy to verify that $D_{5} = i(D_{3}D_{1}-D_{1}D_{3})$.  
    This completes the proof.
\end{proof}

\section{The algebra $\mathcal D(W)$ for Hermite examples}\label{Hermite Sect}

There are only three families $2\times 2$ of irreducible weight matrices of Hermite-type such that their algebra $\mathcal D(W)$ admits a $W$-symmetric second-order differential operator with leading coefficient $F_{2}(x) = I$, see \cite[Theorem 3.1]{DG07}. They are 
\begin{align}  \label{exI}
    W(x) &= e^{-x^{2}} \begin{pmatrix} 1 +a^{2}x^{4} && ax^{2} \\ax^{2} && 1 \end{pmatrix} = \begin{pmatrix} 1 && ax^{2} \\ 0 && 1 \end{pmatrix} \begin{pmatrix} e^{-x^{2}} && 0 \\ 0 && e^{-x^{2}} \end{pmatrix}  \begin{pmatrix} 1 && 0 \\ ax^{2} && 1 \end{pmatrix},\\  \label{exII}
W(x) & = e^{-x^{2}}\begin{pmatrix} 1 + a^{2}x^{2} && ax \\ ax && 1 \end{pmatrix} = \begin{pmatrix} 1 && ax \\ 0 && 1 \end{pmatrix} \begin{pmatrix} e^{-x^{2}} && 0 \\ 0 && e^{-x^{2}} \end{pmatrix}  \begin{pmatrix} 1 && 0 \\ ax && 1 \end{pmatrix},  \intertext{ and}   \label{exIII}
W(x) &= e^{-x^{2}} \begin{pmatrix} e^{2bx} + a^{2}x^{2} && ax \\ax && 1 \end{pmatrix} = \begin{pmatrix} 1 && ax \\ 0 && 1 \end{pmatrix} \begin{pmatrix} e^{-x^{2}+2bx} && 0 \\ 0 && e^{-x^{2}} \end{pmatrix} \begin{pmatrix} 1 && 0 \\ ax && 1 \end{pmatrix}, 
\end{align}
where $x \in \mathbb{R}$, $a,\, b \in \mathbb{R} 
 - \{ 0 \}$. 

The weight $W$ given in \eqref{exII} has been previously studied in \cite{CG06}, and its algebra $\mathcal{D}(W)$  has been determined in  \cite{T11}. 
On the other hand, the weight in \eqref{exI} has not received the same level of attention. 
Almost nothing is known about its algebra, even at the level of conjecture. 
The only available information is that it admits (essentially) a unique second-order differential operator in the algebra $\mathcal D(W)$. 
 The weight in \eqref{exIII}, was fully studied in \cite{BP23}, where we proved that the algebra $\mathcal D(W)$ is a polynomial algebra in the second-order differential operator $D = \partial^{2}I + \partial \begin{psmallmatrix} -2x + 2b && -2abx + 2a \\ 0 && -2x \end{psmallmatrix} + \begin{psmallmatrix} - 2 && 0 \\ 0 && 0  \end{psmallmatrix}.$

\smallskip

This section aims to determine the algebra $\mathcal{D}(W)$ for the weights given in \eqref{exI} and \eqref{exII}, taking full advantage of the fact that the weight $W$ can be obtained as a Darboux transformation of the scalar weight $\widetilde W(x)=e^{-x^2} I$.

\begin{remark}
The weight matrix $W$ in \eqref{exIII} is not a Darboux transformation of a direct sum of scalar weights.   

If fact, let $w_{1},w_{2}$ be scalar weights supported on $\mathbb{R}$. If  $W$ is a Darboux transformation of $\widetilde W(x) = \begin{psmallmatrix} w_{1}(x) && 0 \\0 && w_{2}(x) \end{psmallmatrix}$ then there exist degree-preserving differential operators $\mathcal{V}$ and $\mathcal{N}$ as in the Definition \ref{darb def}.  We observe that the $0$-th order differential operator $\begin{psmallmatrix} 1 && 0 \\ 0 && 0 \end{psmallmatrix}$ belongs to $\mathcal{D}(\widetilde W)$. By Proposition \ref{alg debil}, we have that the differential operator $\mathcal{N}\begin{psmallmatrix} 1 && 0 \\ 0 && 0 \end{psmallmatrix}\mathcal{V} \in \mathcal{D}(W)$. However, its leading coefficient is singular, which implies that it is not a polynomial in $D$.  This gives a contradiction.
\end{remark}



Let us consider a direct sum of scalar Hermite weights, $\widetilde W(x) = e^{-x^{2}}I$. 
  From  \cite{BP23b}, we  know the structure of the algebra of this weight $\widetilde W$ 
\begin{equation} \label{her}
    \mathcal{D}(\widetilde W) = \left \{ \begin{pmatrix}p_{1}(\delta) && p_{2}(\delta) \\ p_{3}(\delta) && p_{4}(\delta) \end{pmatrix}  \biggm| p_{1},\, p_{2},\, p_{3}, \, p_{4} \in \mathbb{C}[x] \right \},
\end{equation} 
where $\delta = \partial^{2} -2 \partial x$ is the  classical  (scalar) Hermite differential operator.
In particular, $\mathcal{D}(\widetilde W)$ is a non-commutative algebra, and every differential operator within it is of even order.

\begin{prop}\label{sim her}
    The leading coefficient of a 
    $\widetilde W$-symmetric differential operator  $D = \displaystyle\sum_{j=0}^{2m}\partial^{j}F_{j} $ is 
    $$F_{2m}(x) = \begin{pmatrix} k_{1} && k_{2} + ik_{3} \\ k_{2}-ik_{3} && k_{4} \end{pmatrix} \text{ for some } k_{1},\,k_{2}, \, k_{3}, \, k_{4} \in \mathbb{R}.$$
\end{prop}
\begin{proof}  
    From the symmetry equations given in  Theorem \ref{equivDsymm}, we obtain that the leading coefficient of $D$ satisfies $F_{2m}(x)\widetilde W(x) = \widetilde W(x)F_{2m}(x)^{\ast}$. Now, the statement follows easily from \eqref{her}.
\end{proof}

\subsection{Example I}
\mbox{}

We consider  the first family mentioned above
$$W(x) = e^{-x^{2}} \begin{pmatrix} 1 + a^{2}x^{4} && ax^{2} \\ ax^{2} && 1 \end{pmatrix}, \qquad a\not=0.
$$

First of all, we will prove here that $W$ is a Darboux transformation of the diagonal weight $\widetilde W(x) = e^{-x^2}$. Consequently, we will provide an explicit expression for a sequence of orthogonal polynomials along with their three-term recurrence relation. Additionally, we will determine the generators of $\mathcal{D}(W)$, establish some algebraic relations among them, and identify the generators of the center $\mathcal{Z}(W)$ within $\mathcal{D}(W)$.

\begin{prop}
    The weight matrix $W(x)$ is a strong Darboux transformation of the direct sum of classical scalar Hermite weights $\widetilde W(x) = e^{-x^{2}}I$.
\end{prop}
\begin{proof} Let us consider  $\delta = \partial^{2} - 2\partial\,  x$  the classical (scalar) Hermite operator. We have that the differential operator 
    $$D = \begin{pmatrix} \delta^{2} - 6\delta +8 +\frac{16}{a^{2}} && 0 \\ 0 && \delta^{2} + 2\delta + \frac{16}{a^{2}} \end{pmatrix} \in \mathcal{D}(\widetilde W) $$
     and it can be factorized as $D = \mathcal{V}\mathcal{N}$, with
    \begin{equation} \label{strong example 1}
        \begin{split}
            \mathcal{V} & =\partial^{2}\begin{pmatrix}0 && 1\\-1 && ax^{2} \end{pmatrix} + \partial \begin{pmatrix} 0 && -4x \\ 0 && 0 \end{pmatrix} + \begin{pmatrix} \frac{4}{a} && -2 \\ 0 && \frac{4}{a} \end{pmatrix}, \\
            \mathcal{N}  = W\mathcal{V}^{\ast}\widetilde{W}^{-1} & = \partial^{2} \begin{pmatrix}ax^{2} && -1 \\1 && 0 \end{pmatrix} + \partial \begin{pmatrix} 4ax && 4x \\ 0 && 0 \end{pmatrix} + \begin{pmatrix} 2a + \frac{4}{a} && 2 \\ 0 && \frac{4}{a} \end{pmatrix}.
        \end{split}
    \end{equation}
   By using integration by parts it is easy to verify that 
  $$\langle P \cdot \mathcal{V}, Q \rangle_{W} = \langle P, Q\cdot \mathcal{N} \rangle_{\widetilde W} \qquad \text{for all } P, \, Q \in \operatorname{Mat}_{2}(\mathbb{C}[x]).$$ 
 Hence, the statement follows by  Proposition \ref{strong1}.
\end{proof}

\begin{remark}  
It is worth noting that the Darboux transformation between $W$ and $\widetilde W$ cannot arise from the factorization of a second-order differential operator $D \in \mathcal{D}(\widetilde W)$. Otherwise, $\mathcal{V}$ and $\mathcal{N}$ would be of order $1$. Then, by Proposition \ref{alg debil}, the algebra $\mathcal{D}(W)$ would have the following second-order operator given by $\mathcal{N}\begin{pmatrix} 1 && 0 \\ 0 && 0 \end{pmatrix} \mathcal{V}$, but its leading coefficient is a singular matrix. We will see later that the leading coefficient of the differential operators of order less than or equal to $2$ in the algebra $\mathcal{D}(W)$ is $\alpha I$ for some $\alpha \in \mathbb{C}$. Thus, it leads us to a contradiction.
\end{remark}

\medskip

The strong Darboux transformer 
gives us an explicit expression of an orthogonal sequence of polynomials for $W$.
\begin{cor}
  A sequence of orthogonal polynomials with respect to $W(x)$ 
 is given by  
$$Q_{n}(x) = \begin{pmatrix} \frac{4}{a}H_{n}(x) && -2H_{n}(x)-4xnH_{n-1}(x)+n(n-1)H_{n-2}(x) \\ -n(n-1)H_{n-2}(x) && \frac{4}{a}H_{n}(x) + ax^{2}n(n-1)H_{n-2}(x)\end{pmatrix},$$
where $H_{n}(x)$ is the sequence of monic orthogonal polynomials for the classical scalar Hermite weight.
  \end{cor}
\begin{proof}
 The monic orthogonal polynomials with respect to    $\widetilde W(x) = e^{-x^2}I$ are given by 
 $P_{n}(x) = \operatorname{diag}(H_{n}(x), H_{n}(x))$.
Since $W(x)$ is a strong Darboux transformation of $\widetilde W(x)$, from Proposition \ref{cons},  we obtain that
 $Q_n(x)=P_{n}(x)\cdot \mathcal{V} $ is a sequence of orthogonal polynomials for $W$.
\end{proof}

\begin{remark}
   We observe that  the leading coefficient of $Q_{n}(x)$ is $\begin{psmallmatrix} \frac{4}{a} && -2-4n \\0 && \frac{4}{a} + an(n-1) \end{psmallmatrix}$. Also, by direct calculation, we obtain the three-term recurrence relation
   $$x Q_n(x)= A_n Q_{n+1}(x)+ B_n Q_n(x)+ C_n Q_{n-1}(x), $$
where
$$ A_n=\begin{pmatrix}  1 &&   \frac{4a}{a^{2}n^{2}+a^{2}n+4} \\  0 &&  \frac{a^{2}n^{2}-a^{2}n+4}{a^{2}n^{2}+a^{2}n+4} \end{pmatrix}, \qquad B_n=0, \quad \text{ and }\quad C_n=\begin{pmatrix}  \frac{n(a^{2}n^{2}+3a^{2}n+2a^{2}+4)}{2(a^{2}n^{2}+a^{2}n+4)} && 0 \\  \frac{2an}{a^{2}n^{2}+a^{2}n+4} &&  \frac{n}{2} \end{pmatrix}.$$
\end{remark}

\

Now, we aim to determine the algebra $\mathcal{D}(W)$. 
From Corollary \ref{parity cor}, we have that every differential operator $D \in \mathcal{D}(W)$ has even order. We use the symmetry equations \eqref{symmeq2} to obtain that the algebra $\mathcal{D}(W)$ contains the  second-order $W$-symmetric differential operator 
$$D_{1}  = \partial^{2} I + \partial \begin{pmatrix} -2x && 4ax \\ 0 &&  -2x \end{pmatrix} + \begin{pmatrix}  0 && 2a  \\ 0 && 4 \end{pmatrix}.$$
Moreover the subspace of $\mathcal{D}(W)$ containing differential operators of order equal to or lower than $2$ is generated by $\{I, D_{1}\}$.
See also \cite{DG07}.

We also obtain,   by Theorem \ref{main}, that $\mathcal{D}(W)$ contains the following fourth-order $W$-symmetric differential operators
\begin{align*}
        D_{2} = \mathcal{N}E_{1}\mathcal{V} & = \partial^{4} \begin{psmallmatrix}  -ax^{2} &&  a^{2}x^{4}-1 \\ -1 && ax^{2} \end{psmallmatrix} + \partial^{3} \begin{psmallmatrix} -8ax &&  8a^{2}x^{3}+8x \\ 0 && 0 \end{psmallmatrix} \\
        & \quad + \partial^{2} \begin{psmallmatrix} -12a-\frac{8}{a} &&  12a^{2}x^{2}-8x^{2}+12 \\ 0 && \frac{8}{a}  \end{psmallmatrix} + \partial \begin{psmallmatrix} \frac{16x}{a} && -16x \\ 0 && -\frac{16x}{a} \end{psmallmatrix} + \begin{psmallmatrix} \frac{8}{a} &&  4+\frac{16}{a^{2}} \\ \frac{16}{a^{2}} &&  -\frac{8}{a}\end{psmallmatrix}, \displaybreak[0]\\
        D_{3} = \mathcal{N}E_{2}\mathcal{V} & = \partial^{4} \begin{psmallmatrix} 1 && -ax^{2} \\  0 &&  0 \end{psmallmatrix} + \partial^{3} \begin{psmallmatrix}  -4x && 4ax^{3} \\ 0  && 0\end{psmallmatrix} + \partial^{2} \begin{psmallmatrix}  -10 && 10ax^{2}-\frac{4}{a} \\ -\frac{4}{a} &&  4x^{2} \end{psmallmatrix}
        + \partial \begin{psmallmatrix}  0 && \frac{16x}{a} \\ 0  && 0 \end{psmallmatrix} + \begin{psmallmatrix}  0 &&  \frac{8}{a} \\ 0 && \frac{16}{a^{2}}\end{psmallmatrix},  \\
        D_{4} = \mathcal{N}E_{3}\mathcal{V} & = \partial^{4} \begin{psmallmatrix} 0 && ax^{2} \\ 0 && 1 \end{psmallmatrix} + \partial^{3} \begin{psmallmatrix}  0 &&  -4ax^{3}+8ax \\  0 && -4x \end{psmallmatrix} + \partial^{2} \begin{psmallmatrix}  4x^{2} &&  -26ax^{2}+12a+\frac{4}{a} \\  \frac{4}{a} &&  -2 \end{psmallmatrix} \\
        & \quad + \partial \begin{psmallmatrix} 16x && -\frac{16x(2a^{2}+1)}{a} \\  0 && 0\end{psmallmatrix} + \begin{psmallmatrix} \frac{8(a^{2}+2)}{a^{2}} && -\frac{4(a^{2}+2)}{a} \\ 0 && 0 \end{psmallmatrix}, \\
        D_{5} = \mathcal{N}E_{4}\mathcal{V} & = \partial^{4}i\begin{psmallmatrix} -ax^{2} && a^{2}x^{4}+1 \\ -1 && ax^{2} \end{psmallmatrix} + \partial^{3} i \begin{psmallmatrix} -8ax && 8a^{2}x^{3}-8x \\ 0 && 0 \end{psmallmatrix} \\
        & \quad + \partial^{2} i \begin{psmallmatrix} -12a && 12a^{2}x^{2}+24x^{2}-12 \\ 0 && 0\end{psmallmatrix} + \partial i \begin{psmallmatrix} -\frac{16x}{a} && 48x \\ 0 && \frac{16x}{a} \end{psmallmatrix} + i \begin{psmallmatrix}  -\frac{8}{a} && 12+\frac{16}{a^{2}} \\ -\frac{16}{a^{2}} && \frac{8}{a}\end{psmallmatrix},  
\end{align*}
where $E_{1} = \begin{pmatrix} 0 && 1 \\ 1 && 0 \end{pmatrix}$, $E_{2} = \begin{pmatrix} 0 && 0 \\ 0 && 1 \end{pmatrix}$, $E_{3} = \begin{pmatrix}1 && 0 \\ 0 && 0  \end{pmatrix}$ and $E_{4} = \begin{pmatrix}0 && i \\ -i && 0 \end{pmatrix}$ which are $\widetilde W$-symmetric differential operators.

\begin{prop}
The algebra $\mathcal{D}(W)$ is generated by $\{ I, D_{1}, D_{2}, D_{3} \}$.
\end{prop}

\begin{proof} 
For any $W$-symmetric differential operator $D = \sum_{j=0}^{2n}\partial^{j}F_{j} $
in $\mathcal{D}(W)$ we have that $\mathcal{V}D\mathcal{N} \in \mathcal{S}(\widetilde W)$, where $\mathcal{V}$ and $\mathcal{N}$ are the differential operators in \eqref{strong example 1}. The leading coefficient of $\mathcal{V}$ is $\begin{psmallmatrix}0 && 1\\-1 && ax^{2} \end{psmallmatrix}$, 
of $\mathcal{N}$ is $\begin{psmallmatrix} ax^{2} && -1\\1 && 0 \end{psmallmatrix}$ and the leading coefficient of $\mathcal{V}D\mathcal{N}$ is  
$\begin{psmallmatrix} k_{1} && k_{2} + i k_{3} \\ k_{2} - ik_{3} && k_{4} \end{psmallmatrix}$, 
for some $k_{1},k_{2},k_{3},k_{4} \in \mathbb{R}$.  Thus we obtain that the leading coefficient of $D$ is
\begin{equation}\label{F2n}
\begin{split}
F_{2n} &= \begin{pmatrix} ax^{2} && -1\\1 && 0 \end{pmatrix}\begin{pmatrix} k_{1} && k_{2} + i k_{3} \\ k_{2} - ik_{3} && k_{4} \end{pmatrix}\begin{pmatrix} 0 && 1 \\-1 && ax^{2} \end{pmatrix}\\
&=\begin{pmatrix} -ik_{3}ax^{2}-k_{2}ax^{2}+k_{4} && ik_{3}a^{2}x^{4}+k_{2}a^{2}x^{4}+k_{1}ax^{2}-k_{4}ax^{2}-k_{2}+ik_{3} \\  -k_{2}-ik_{3} && ik_{3}ax^{2}+k_{2}ax^{2}+k_{1} \end{pmatrix}.
\end{split}
\end{equation}

Now, we  prove that $\mathcal{D}(W)$ is generated by $\{I,D_{1},D_{2},D_{3},D_{4},D_{5} \}$. Let $D=\sum_{j=0}^{2n}\partial^{j}F_{j} \in \mathcal{S}(W)$.  
For  $n \leq 1$ the statement is true. Suppose that $n > 1$. 
 The $W$-symmetric differential operator $$\widetilde{E} = D_{1}^{n-2}E + ED_{1}^{n-2}, \quad \text{  with } \quad E = \frac{1}{2}(k_{2}D_{2} + k_{4}D_{3} + k_{1}D_{4}+k_{3}D_{5}),$$ is of order $2n$ and its  leading coefficient is $F_{2n}$. Thus $D - \widetilde{E}$ is a $W$-symmetric differential operator of order $< 2n$.  
Hence, by the inductive hypothesis, we have that  $D$ is generated by $\{I,D_{1},D_{2},D_{3},D_{4}, D_5 \}$. 
Finally, by computing  the eigenvalues of the differential operators, we get 
\begin{align*}
        \Lambda_{n}(D_{4}) & = \Lambda_{n}(D_{1})^{2} -6\Lambda_{n}(D_{1})-\Lambda_{n}(D_{3})+\tfrac{8(a^{2}+2)}{a^{2}}I, \\
        \Lambda_{n}(D_{5}) &= -\frac{i}{4}\big(\Lambda_{n}(D_{2})\Lambda_{n}(D_{1}) - \Lambda_{n}(D_{1})\Lambda_{n}(D_{2})\big).
\end{align*}
Hence, from Proposition \ref{prop2.8-GT}, $D_{4}$ and $D_{5}$ are generated by $\{I, D_{1}, D_{2}, D_{3} \}$.
\end{proof}

The eigenvalues of $D_{1}$, $D_{2}$, and $D_{3}$ are

\begin{equation}
    \begin{split}
        \Lambda_{n}(D_{1}) & = \begin{pmatrix} -2n && 4an+2a \\ 0 && -2n + 4 \end{pmatrix}, \\
        \Lambda_{n}(D_{2}) & = \begin{pmatrix} \frac{8(2n+1)}{a} && \frac{(a^{2}n^{2}+2a^{2}n-4)(a^{2}n^{2}-a^{2}-4)}{a^{2}} \\ \frac{16}{a^{2}} && -\frac{8(2n+1)}{a} \end{pmatrix}, \\
        \Lambda_{n}(D_{3}) & = \begin{pmatrix} 0 &&  \frac{(4n+2)(a^{2}n^{2}-a^{2}n+4)}{a} \\   0 && \frac{(4a^{2}n^{2}-4a^{2}n+16)}{a^{2}} \end{pmatrix}.
    \end{split}
\end{equation}

\begin{remark}
    The differential operators $D_{1},D_{2}$ and $D_{3}$ satisfy the following relations
\begin{equation*}
    \begin{split}
            D_{1}D_{3} & = D_{3}D_{1} = \tfrac{1}{16}\left(D_{1}^{2} + 2D_{1} + \tfrac{16}{a^2}\right )\left(D_{1}^{2}-6D_{1}+8+\tfrac{16}{a^2}\right )-\tfrac{1}{16}D_{2}^{2} + 3D_{3}, \\
            D_{2}D_{3} & = -\tfrac{1}{8}\left(D_{1}D_{2}-D_{2}D_{1}-4D_{2}\right )\left (D_{1}^{2}-6D_{1}+\tfrac{8(a^{2}+2)}{a^{2}}I\right ), \\
            D_{3}D_{2} & = D_{1}^{2}D_{2}-6D_{1}D_{2} + \tfrac{8(a^{2}+2)}{a^{2}}D_{2}+\tfrac{1}{8}\left(D_{1}^{2}-6D_{1}+\tfrac{8(a^{2}+2)}{a^{2}}I\right )\left(D_{1}D_{2}-D_{2}D_{1}-4D_{2}\right).
        \end{split}
    \end{equation*}
\end{remark}

\

We want to determine the center $\mathcal{Z}(W)$ of the algebra $\mathcal{D}(W)$. 
\begin{prop}\label{center}
    If  $D \in \mathcal{Z}(W)$, then the leading coefficient of $D$ is $\alpha I$ for some $\alpha \in \mathbb{C}$.
\end{prop}
\begin{proof}
    Given a differential operator $D \in \mathcal{Z}(W)$, it follows that its adjoint $D^{\dagger}$ belongs to $\mathcal{Z}(W)$ because ${\dagger}$ is an involution. Thus, it is easy to prove that $\mathcal{Z}(W) = \left ( \mathcal{Z}(W)\cap \mathcal{S}(W) \right ) \oplus i\left(\mathcal{Z}(W)\cap\mathcal{S}(W) \right)$. Then, it is sufficient to prove the statement for the symmetric operators in the center of $\mathcal{D}(W)$. 
    
    If $D = \sum_{j=0}^{2m}\partial^{j}F_{j} \in \mathcal{Z}(W)\cap\mathcal{S}(W)$ then $D$ commutes with the generators of the algebra $\mathcal{D}(W)$. In particular,
    \begin{equation} \label{pip}
        DD_{2} = D_{2}D \quad \text{ and } \quad DD_{3} = D_{3}D.
    \end{equation}
    By comparing leading coefficients in \eqref{pip}, it follows that $F_{2m}$, the leading coefficient of $D$, is of the form $p(x)I$ for some polynomial $p(x)$. On the other hand, 
    we have that $F_{2m}$ is as in \eqref{F2n}. Thus, we conclude that $F_{2m}(x) = k_{1}I$ for some $k_{1} \in \mathbb{R}$.
\end{proof}

Let us define
\begin{equation*}
    \begin{split}
        Z_{1} & = D_{1}^{3} - 3D_{1}^{2} - \tfrac{(4a^{2}-48)}{a^{2}}D_{1} - 12 D_{3}, \\
        Z_{2} & = D_{2}^{2}, \\
        Z_{3} & = D_{1}^{5} -20D_{3}^{2} - 80D_{3}D_{1} - \tfrac{20(a^{2}-8)}{3a^{2}}D_{3} - 40 D_{1}^{2} + \tfrac{32(a^{4}+40a^{2}+120)}{3a^{2}}D_{1}.
    \end{split}
\end{equation*}
We have that the eigenvalues of $Z_{1},Z_{2}$, and $Z_{3}$ are
\begin{equation*}
    \begin{split}
    \Lambda_{n}(Z_{1}) & = -4 \, {n \left( 2n^{2}+3n-2+\tfrac{24}{a^2}
 \right) }I, \\
    \Lambda_{n}(Z_{2}) & = 16 \,  \left( {n}^{2}+3n+2+\tfrac{4}{a^2} \right) 
 \left( {n}^{2}-n+\tfrac{4}{a^2} \right) I, \\
 \Lambda_{n}(Z_{3}) & = -\tfrac {32}{3}\,n \left( 3n^4-5{n}^{2}
+15n+2+\tfrac{40}{a^2}\,{n}^{2}+\tfrac{80} {{a}^{2}}+\tfrac{240} {a^{4}}\right) I.
    \end{split}
\end{equation*}
For $D \in \mathcal{D}(W)$ we have $\Lambda_{n}(DZ_{i}) = \Lambda_{n}(D)\Lambda_{n}(Z_{i}) = \Lambda_{n}(Z_{i})\Lambda_{n}(D) = \Lambda_{n}(Z_{i}D)$ for all $n \geq 0$. 
By Proposition \ref{prop2.8-GT}, it follows that $DZ_{i} = Z_{i}D$. Consequently, $Z_{1}, Z_{2}, Z_{3}$ belong to $\mathcal{Z}(W)$. 

\begin{prop}
    The center $\mathcal{Z}(W)$ of $\mathcal{D}(W)$ is generated by $\{I,Z_{1},Z_{2},Z_{3} \}$.
\end{prop}
\begin{proof} 
Let us prove that there are no operators of order 2 or 4 in $\mathcal{Z}(W)$. 
If such an operator existed, its leading coefficient would be $\alpha I$, for some $\alpha \in \mathbb C$.  Then 
it would be a linear combination of $D_1^2$, $D_1$, and $I$. 
But, $\beta_1 D_{1}^{2} + \beta_2 D_{1} + \mu I$ does not commute with $D_{2}$ for any $\beta_1, \beta_2 \neq 0$. 

    Recall that every differential operator $D \in \mathcal{D}(W)$ has even order and let $D = \sum_{j=0}^{2m}\partial^{j}F_{j} \in \mathcal{Z}(W)$. By Proposition \ref{center} we have that $F_{2m} = \alpha I$ for some $\alpha \in \mathbb{C}$.  We proceed by induction on $m$. For $m = 3$, we have that $D$ is of order $6$ and $F_{6} = \alpha I$ for some $\alpha \in \mathbb{C}$. Then, $D - \alpha Z_{1}$ belongs to $\mathcal{Z}(W)$ and it is of order less than $6$. Thus, $D = \alpha Z_{1} + kI$ for some $k\in \mathbb{C}$. 
    
    For $m > 3$, we have that $Z_{1},Z_{2}$ and $Z_{3}$ are differential operators of order $6$, $8$ and $10$ respectively. Thus, we can find integers $s,t,r \in \mathbb{N}_{0}$ such that $6s + 8t + 10r = 2m$. Then $B = D - \alpha Z_{1}^{s}Z_{2}^{t}Z_{3}^{r}\in \mathcal{Z}(W)$ is of order less than $2m$. Therefore, by the inductive hypothesis, the differential operator $B$ is generated by $\{ I, Z_{1}, Z_{2}, Z_{3} \}$ and hence also the differential operator $D$.
\end{proof}

\subsection{Example II}
\mbox{}

Now we consider  the weight matrix
$$W(x) = e^{-x^{2}}\begin{pmatrix} 1 + a^{2}x^{2} && ax \\ax && 1 \end{pmatrix}, \quad a \in \mathbb{R} - \{0\}, \, x\in \mathbb{R}.$$

As we mentioned before, the algebra $\mathcal{D}(W)$ associated with this weight has already been studied by J. Tirao in \cite{T11},  a comprehensive and technical paper. Here, we provide an alternative and shorter method to compute this algebra, as a consequence of the fact that the weight $W$ is a Darboux transformation of $\widetilde W$.


\begin{prop} The weight $W(x)$ is a strong Darboux transformation of $\widetilde W(x) = e^{-x^{2}}I$.
\end{prop}
\begin{proof}
    We take the differential operator 
    $$D = \begin{pmatrix} -\delta + 2+\frac{4}{a^{2}} && 0 \\ 0 && -\delta + \frac{4}{a^{2}} \end{pmatrix} \in \mathcal{D}(\widetilde W)$$
     where $\delta = \partial^{2} + \partial (-2x)$ is the second-order differential operator of the classical scalar Hermite weight, and we factorize it as $D = \mathcal{V}\mathcal{N}$, with
     \begin{equation} \label{strong example 2}
        \begin{split}
            \mathcal{V} & = \partial \begin{pmatrix} 0 && 1 \\ 1 && -ax \end{pmatrix} - \frac{2}{a}I, \\
            \mathcal{N}  = W\mathcal{V}^{\ast}\widetilde{W}^{-1} & =   \partial \begin{pmatrix} -ax && -1 \\ -1 && 0 \end{pmatrix} + \begin{pmatrix} -a - \frac{2}{a} && 0 \\ 0 && -\frac{2}{a} \end{pmatrix}.
        \end{split}
    \end{equation}
    
  By integration by parts, we easily obtain that $\langle P \cdot \mathcal{V}, Q \rangle_{W} = \langle P, Q\cdot \mathcal{N} \rangle_{\widetilde W}$,  for all $P, \, Q \in \operatorname{Mat}_{2}(\mathbb{C}[x])$.  
  Thus, by Proposition \ref{strong1}, the statement holds. 
 \end{proof}
 
 From Theorem \ref{main} we can map $\widetilde W$-symmetric operators into $W$-symmetric operators. Thus we have the following (second-order) $W$-symmetric differential operators in the algebra $\mathcal D(W)$.
\begin{equation*}
    \begin{split}
        D_{1}& = \mathcal{N}E_{1}\mathcal{V}  =  \partial^{2} \begin{pmatrix} -ax &&  a^{2}x^{2}-1 \\  -1 &&  ax \end{pmatrix} + \partial \begin{pmatrix} -2a && 2a^{2}x+4x \\ 0 &&  0 \end{pmatrix} + \begin{pmatrix} 0 &&  2+\frac{4}{a^{2}} \\ \frac{4}{a^{2}} &&  0 \end{pmatrix}, \\
        D_{2} &= \mathcal{N}E_{2}\mathcal{V} =  \partial^{2} \begin{pmatrix} -1 &&  ax \\ 0 &&  0 \end{pmatrix} + \partial \begin{pmatrix}  0 &&  \frac{2}{a} \\  -\frac{2}{a} &&  2x\end{pmatrix} + \begin{pmatrix} 0 &&  0 \\ 0 &&  \frac{4}{a^{2}}\end{pmatrix}, \\
        D_{3} & = \mathcal{N}E_{3}\mathcal{V} = \partial^{2} \begin{pmatrix} 0 &&  -ax \\ 0 && -1 \end{pmatrix} + \partial \begin{pmatrix}  2x && -2a - \frac{2}{a} \\ \frac{2}{a} && 0\end{pmatrix} + \begin{pmatrix}  2 + \frac{4}{a^{2}} &&  0 \\ 0  && 0\end{pmatrix}, \\
        D_{4} & = \mathcal{N}E_{4}\mathcal{V} = \partial^{2}i \begin{pmatrix}  -ax && a^{2}x^{2}+1 \\  -1 && ax\end{pmatrix} + \partial i \begin{pmatrix} -2a-\frac{4}{a} &&  2a^{2}x+4x \\  0 && \frac{4}{a}\end{pmatrix} + i \begin{pmatrix} 0 && 2+\frac{4}{a^{2}} \\  -\frac{4}{a^{2}} &&  0 \end{pmatrix}, 
    \end{split}
\end{equation*}
where $E_{1}= \begin{psmallmatrix} 0 && 1 \\ 1 && 0 \end{psmallmatrix}$, $E_{2} = \begin{psmallmatrix} 0 && 0 \\ 0 && 1 \end{psmallmatrix}$, $E_{3} = \begin{psmallmatrix} 1 && 0 \\ 0 && 0 \end{psmallmatrix}$ and $E_{4} = i\begin{psmallmatrix} 0 && 1 \\ -1 && 0 \end{psmallmatrix}$, and $\mathcal{V}$, $\mathcal{N}$ are defined in \eqref{strong example 2}.
We also introduce the differential operator $D_5  = -\mathcal{N}\mathcal{V}$, whose leading coefficient is a scalar matrix
$$ D_{5}  = -(D_{2} + D_{3}) = \partial^{2} I + \partial \begin{pmatrix}  -2x &&  2a \\ 0 &&  -2x \end{pmatrix} + \begin{pmatrix} -2-\frac{4}{a^{2}} && 0 \\ 0 && -\frac{4}{a^{2}} \end{pmatrix}. $$

\begin{thm}
    The algebra $\mathcal{D}(W)$ is generated by $\{ I, D_{1}, D_{5} \}$.
\end{thm}
\begin{proof}
The weight $W$ is a strong Darboux transformation of $e^{-x^{2}}I$, thus by Corollary \ref{parity cor}, we have that all differential operators in $\mathcal{D}(W)$ are of even order.
If  $D = \sum_{j=0}^{2n}\partial^{j}F_{j} \in \mathcal{S}(W)$ then $\mathcal{V}D\mathcal{N} \in \mathcal{S}(\widetilde W)$, and thus, from Proposition \ref{sim her}, the leading coefficient of $D$ is of the form
\begin{equation*} \label{F2n bis}
\begin{split}
F_{2n}&= \begin{pmatrix} ax && 1 \\ 1 && 0 \end{pmatrix}\begin{pmatrix} k_{1} && k_{2}+ik_{3} \\ k_{2}-ik_{3} && k_{4} \end{pmatrix}\begin{pmatrix} 0 && -1 \\ -1 && ax \end{pmatrix}
\\ &= \begin{pmatrix}  -k_{2}ax-ik_{3}ax-k_{4} && k_{2}a^{2}x^{2}+ik_{3}a^{2}x^{2}-k_{1}ax+k_{4}ax-k_{2}+ik_{3} \\ -k_{2} - ik_{3} &&  k_{2}ax+ik_{3}ax-k_{1} \end{pmatrix}.  
\end{split}
\end{equation*}

We observe  that $F_{2n}$ is  a linear combination of the leading coefficients of $D_{1}, D_{2}, D_{3}$ and $D_{4}$.  Additionally we can  use the operator $D_{5}$ 
to increase the order of any differential operator without modifying its leading coefficient. Therefore, we can prove by induction that $\mathcal{D}(W)$ is generated by $\{ D_{1},D_{2},D_{3},D_{4},D_{5} \}$. Finally, by observing that 
$$ D_{3}  = \tfrac{1}{4}(D_{5}^{2} - D_{1}^{2}) - \frac{1}{2}D_{5}, \qquad 
        D_{2}  = -D_{3}- D_{5}, \qquad 
        D_{4}  = \frac{i}{2}(D_{1}D_{5}-D_{5}D_{1}).$$
we obtain that the algebra $\mathcal{D}(W)$ is generated by $\{I,D_{1},D_{5} \}$.
\end{proof}

From Proposition \ref{strong1} we have that $Q_{n}(x) = H_{n}(x)\cdot \mathcal{V}$ is an orthogonal sequence of polynomials for $W$. Explicitly, we have
\begin{equation} \label{II-Hermpoly}Q_{n}(x) = \begin{pmatrix} -\frac{2}{a}H_{n} && nH_{n-1} \\ nH_{n-1} && -\frac{2}{a}H_{n} - anxH_{n-1} \end{pmatrix},
\end{equation}
where $H_{n}(x)$ is the sequence of monic orthogonal polynomials for the classical scalar Hermite weight. Observe that the leading coefficient of $Q_{n}(x)$ is $\begin{psmallmatrix} -\frac{2}{a} && 0 \\ 0 && -\frac{2}{a} - an \end{psmallmatrix}.$

\section{Beyond the Matrix Bochner Problem: New weight matrices}
Throughout this paper, we have illustrated some profound relationships that exist between the algebras of differential operators associated with weight matrices that are Darboux transformations of each other.
In this section, we introduce two new families of irreducible Hermite-type weight matrices, of size $2\times 2$. Both examples can be obtained as Darboux transformations of a direct sum of classical weights $\widetilde W(x) = e^{-x^{2}}I$. However, they are not solutions to the Matrix Bochner Problem, meaning that their algebra $\mathcal D(W)$ does not contain second-order differential operators. 

The orthogonal polynomials associated with these weights $W$ are simultaneous eigenfunctions of differential operators of order greater than 2.
To the best of our knowledge, they are the first examples of this situation in the matrix case.
In other words, we present examples of irreducible weight matrices where the algebra $\mathcal{D}(W)$ is non-trivial, but it does not contain second-order differential operators.

\medskip

Let us consider the  following Hermite-type weight matrix
\begin{equation}\label{WnoBochner1}
    W(x) = e^{-x^{2}}\begin{pmatrix} a^{4}x^{4}+3a^{2}x^{2}+1 && a^{3}x^{3}+2ax \\ a^{3}x^{3}+2ax && a^{2}x^{2}+1\end{pmatrix}, \quad a \in \mathbb{R}-\{0\}.
\end{equation}

By using the symmetry equations given in \eqref{symmeq2}, we can see that there are no  $W$-symmetric second-order differential operators. 
Therefore, the weight $W$ is not a solution to the Matrix Bochner Problem.  
Nevertheless, $\mathcal D(W)$ contains 
the following fourth-order differential operator

\begin{equation*}
    \begin{split}
        D & = \partial^{4}I + \partial^{3} \begin{pmatrix} 4a^{2}x-4x && -4a^{3}x^{2}+4a \\ 4a &&  -4a^{2}x-4x \end{pmatrix} + \partial^{2} \begin{pmatrix} 12a^{2}+4x^{2}-\frac{8}{a^{2}} && -12a^{3}x-24ax \\ 0 && 4x^{2}-12-\frac{8}{a^{2}} \end{pmatrix} \\
        & \quad + \partial \begin{pmatrix}  12x+\frac{16x}{a^{2}} &&  -24a-\frac{32}{a} \\ 0 && 4x+\frac{16x}{a^{2}}\end{pmatrix} + \begin{pmatrix}  4+\frac{16}{a^{2}} && 0 \\ 0 &&  0\end{pmatrix}.
    \end{split}
\end{equation*}

The weight $W$ can be obtained  through a Darboux transformation of $\widetilde W(x) = e^{-x^{2}}I$, a direct sum of classical scalar Hermite weights.
Let   $\delta = \partial^{2} -2 \partial x$ be the classical Hermite second-order differential operator.

The following eighth-order $\widetilde{W}$-symmetric differential operator 
$$\mathcal D = \begin{psmallmatrix}{\frac { \left( {a}^{ 2}\delta-4 \right)  \left( {a}^{ 6}{\delta}^{3}-6 \left(
{a}^{6}+2\,{a}^{4} \right) {\delta}^{2}+ 8\left( {a}^{6}+6\,{a}^{4}+6
\,{a}^{2} \right) \delta-16(3a^{4}+6a^{2}+4) \right) }{{a}^{8}}} & 0 \\ 0 &  {\frac { \left( {a}^{6}{\delta}^{3}-12\,{a}^{4}{\delta}^{2}+ \left( -4\,{a}^
{6}+48\,{a}^{2} \right) \delta - 64 \right)  \left( {a}^{2}\delta-2\,{a}^{2}-4
 \right) }{{a}^{8}}}\end{psmallmatrix}$$
belongs to the algebra $\mathcal{D}(\widetilde W)$ and it can be factored as $\mathcal D = \mathcal{V} \mathcal N $, where $\mathcal N=W\mathcal{V}^{\ast}\widetilde W^{-1}$ and 
\begin{equation*}
    \begin{split}
        \mathcal{V} & = \partial^{4} \begin{pmatrix}1 && -ax \\  -ax &&  a^{2}x^{2}+1\end{pmatrix} + \partial^{3}\begin{pmatrix} -2x && 2ax^{2}-\frac{4}{a} \\  2ax^{2}-\frac{4}{a} &&  -2a^{2}x^{3}+2x \end{pmatrix} \\ & \quad + \partial^{2}\begin{pmatrix}-6 && \frac{6(a^{2}+2)}{a}x \\ \frac{6(a^{2}+2)}{a}x &&  -6(a^{2}+2)x^{2}-6 \end{pmatrix} + \partial \begin{pmatrix} -\frac{8}{a^{2}}x &&  \frac{(12a^{2}+16)}{a^{3}} \\ \frac{(12a^{2}+16)}{a^{3}} &&  -\frac{(12a^{2}+24)}{a^{2}}x \end{pmatrix} + \begin{pmatrix} -\frac{16}{a^{4}} &&   0 \\  0 &&  -\frac{16}{a^{4}} \end{pmatrix}.
    \end{split}
\end{equation*}

It is a matter of integration by parts to verify that $$\langle P \cdot \mathcal{V}, Q \rangle_{W} = \langle P, Q\cdot \mathcal{N} \rangle_{\widetilde W},\quad \text{ for all $P, \, Q \in \operatorname{Mat}_{2}(\mathbb{C}[x])$}.$$  
In this way, by Proposition \ref{strong1}, we obtain that the weight  $W$, given in \eqref{WnoBochner1}, is a Darboux transformation of $\widetilde{W}(x) = e^{-x^{2}}I$. As a consequence, $Q_{n}(x) = H_{n}(x)\cdot \mathcal{V}$ is a sequence of orthogonal polynomials for $W(x)$.

\

Similarly, we also obtain that the Hermite-type weight matrix
\begin{equation}\label{WnoBochner2}
    W(x) =e^{-x^{2}} \begin{pmatrix} (a x^{2} + b x)^{2}+1 &&  -x(a x+b) \\  -x(a x+b) && 1 \end{pmatrix} , \qquad a,b \in \mathbb{R}-\{0\},
\end{equation}
is a strong Darboux transformation of $\widetilde W(x) = e^{-x^{2}}I$ by factoring the eighth-order $\widetilde{W}$-symmetric differential operator 
$$\mathcal{D}=\begin{psmallmatrix}{\delta}^{4}-2{\delta}^{3}+{\frac {4 \left( -{a}^{2}+8 \right) {\delta}^{2}}{{a}^
{2}}}+{\frac { 8\left( {a}^{4}-20\,{a}^{2}-8\,{b}^{2} \right) \delta}{{a}
^{4}}}+{\frac {64(3\,{a}^{2}+2\,{b}^{2}+4)}{{a}^{4}}}
 & -{\frac {24\,{a}^{2}{\delta}^{2}-48\,{a}^{2}\delta-128}{{a}^{3}}} \\ -{\frac {24\,{a}^{2}{\delta}^{2}-48\,{a}^{2}\delta-128}{{a}^{3}}} & {\delta}^{4}-6\,{\delta}^{3}+{\frac {8 \left( {a}^{2}+4 \right) {\delta}^{2}}{{a}^{
2}}}+{\frac { 32\left( {a}^{2}-2\,{b}^{2} \right) \delta}{{a}^{4}}}+\frac{256}
{{a}^{4}}
 \end{psmallmatrix} \in \mathcal{D}(\widetilde W)$$
with strong Darboux transformer from $\widetilde W$ to $W$ given by 
\begin{equation*}
    \begin{split}
        \mathcal{V} & = \partial^{4}\begin{pmatrix} 1 &&  x(a x+b) \\  0 && 1 \end{pmatrix} + \partial^{3}\begin{pmatrix} -2x &&  -2 a x^{3}-2 b x^{2} \\ 0 &&  -6x  \end{pmatrix} + \partial^{2} \begin{pmatrix}  -6 &&  -6 a x^{2}-6 b x+ \frac{8}{a} \\  -\frac{8}{a} &&  4x^{2}-\frac{8 b x}{a}-6\end{pmatrix} \\
        & \quad + \partial \begin{pmatrix}  0 && -\frac{(24a x+8 b)}{a^{2}} \\ -\frac{(-8a x+8 b)}{a^{2}} &&  -\frac{(8 b^{2}-12 a^{2})x}{a^{2}}\end{pmatrix} + \begin{pmatrix} -\frac{16}{a^{2}} &&  -\frac{16}{a} \\  0 && -\frac{16}{a^{2}} \end{pmatrix}.
    \end{split}
\end{equation*}
A sequence of orthogonal polynomials for $W(x)$ is given by $Q_{n}(x) = H_{n}(x)\cdot \mathcal{V}$.

\smallskip

 The algebra $\mathcal{D}(W)$ contains no differential operator of order $2$ or $4$. However, it  contains the sixth-order $W$-symmetric differential operator
\begin{equation*}
    \begin{split}
        D & = \partial^{6}I + \partial^{5} \begin{pmatrix}  -6x &&  -12ax-6b \\ 0 && -6x\end{pmatrix} + \partial^{4} \begin{pmatrix}  12x^{2}-24 &&  42ax^{2}+12bx-30a \\  0 &&  12x^{2}-6\end{pmatrix} \\
        & \quad + \partial^{3} \begin{pmatrix}  -4x(2x^{2}-21) && -24ax^{3}+24bx^{2}+168ax+24b \\ 0  && -4x(2x^{2}-3)) \end{pmatrix} \\
        & \quad + \partial^{2} \begin{pmatrix} -\frac{(48a^{2}x^{2}-24abx-108a^{2}-48)}{a^{2}} && -\frac{(108a^{2}x^{2}-24b^{2}x^{2}-72abx-72a^{2}+24)}{a} \\  -\frac{24}{a} && -\frac{(24a^{2}x^{2}+24abx-48)}{a^{2}} \end{pmatrix} \\
        & \quad + \partial \begin{pmatrix}  -\frac{(48a^{2}x-48ab+96x)}{a^{2}} &&  -\frac{(72a^{3}x-48ab^{2}x+96ax+96b)}{a^{2}} \\ 0 && -\frac{96x}{a^{2}} \end{pmatrix} + \begin{pmatrix}   -\frac{96}{a^{2}} &&  -\frac{48}{a} \\  0 &&  0\end{pmatrix}.
    \end{split}
\end{equation*}

\end{document}